\documentclass{amsart}
\usepackage{amssymb, bbm, calligra, mathrsfs}
\usepackage[all]{xy}

\newcommand{\rdg}{\hfill$\Box$}

\def\calli#1{\textup{\!\textcalligra{#1}}}

\let\x\times
\let\ho\simeq

\let\si\varsigma

\let\d\partial
\let\then\Rightarrow
\let\wt\widetilde

\def\tam{\tau_{\leq m}}

\def\K{\mathbb K}
\def\1{^{-1}}

\def\c{{\boldsymbol C}}

\def\ad{{\sf ad}}

\DeclareMathOperator\id{id} \DeclareMathOperator\Ob{Ob}
 \DeclareMathOperator\Aut{Aut}

\def\brk#1{\langle#1\rangle}

\newtheorem{Pro}{Proposition}[subsection]
\newtheorem{Le}[Pro]{Lemma}
\newtheorem{The}[Pro]{Theorem}
\newtheorem{Co}[Pro]{Corollary}
\theoremstyle{definition}
\newtheorem{De}[Pro]{Definition}
\theoremstyle{remark}

\newtheorem{Rem}[Pro]{Remark}

\newcommand{\Lem}[1]{\begin{Le}#1\end{Le}}
\newcommand{\Cor}[1]{\begin{Co}#1\end{Co}}

\newcommand{\Z}{{\mathbb Z}}

\def\hom{{\calli{Hom}}}
\let\aa\alpha
\let\ee\epsilon

\def\ta{{\cal T}}
\let\ho\simeq

\def\c{{\cal C}}

\let\x\times

\let\t\otimes

\def\Tor{\mathop{\sf Tor}\nolimits}
\def\der{\mathop{\sf Der}\nolimits}
\def\ker{\mathop{\sf Ker}\nolimits}
\def\Cok{\mathop{\sf Coker}\nolimits}
\def\H{{\mathsf H}}

\def\Sh{\mathop{\sf Shukla}\nolimits}
\def\HML{\mathop{\sf HML}\nolimits}

\let\cal\mathscr

\let\x\times

\let\os\oplus

\let\then\Rightarrow

\let\then\Rightarrow

\def\brk#1{\left\langle#1\right\rangle}

\def\hog#1{[\![#1]\!]}

\def\id{{\mathrm{Id}}}

\def\1{^{-1}}
\def\mod{\textrm{-}\mathbf{mod}}

\DeclareMathOperator\Ch{Ch} 
\DeclareMathOperator\Der{\sf Der}
\DeclareMathOperator\Ext{\sf Ext} 

\DeclareMathOperator\Hom{\sf Hom} 
 
\DeclareMathOperator\Ker{\sf Ker} 

\def\A{{\mathbb A}}

\def\F{{\mathbb F}}

\def\Z{{\mathbb Z}}

\def\c{{\mathscr C}}

\begin{document}
\bigskip

\centerline{\bf { SHUKLA COHOMOLOGY AND ADDITIVE TRACK THEORIES} }
\bigskip

\bigskip
\bigskip
\bigskip

\centerline{\bf { Hans-Joachim BAUES and  Teimuraz PIRASHVILI}}
\bigskip
\bigskip
\bigskip

\bigskip
\bigskip
\section{Introduction}
 It is well known that the Hochschild cohomology for associative algebras has
good properties only for algebras which  are projective modules over the ground ring. For
general algebras behavior of Hochschild cohomology is more pathological, for
example there is no long cohomological exact  sequence corresponding to a short
exact sequence for coefficients, etc. In early 60-s Shukla \cite{sh} developed a
cohomology theory for associative algebras with nicer properties than Hochschild theory. Quillen in \cite{Q}
indicated that   the Shukla cohomology fits in his  general framework of
homotopical algebra. The approach of Quillen is based on simplicial methods,
which are usually quite hard to deal with. The aim of this work is to give the foundation of
Shukla cohomology based on chain algebras. We also give an application to the
problem of strengthening additive track theories, which is based on the
comparison homomorphism between Shukla and Mac Lane cohomology theories
\cite{PW}.
We believe that our approach is much
 simpler than one used in \cite{Q} or \cite{sh}.

Let us recall that  a track category is a category enriched in groupoids. A
track category $\ta$ is called \emph{abelian} if for any arrow $f$ the group of automorphisms
of $f$ is abelian. A track theory is an abelian track category with finite lax
products. If it admits strong products then it is called strong track theory.
The main result of \cite{BJP}  asserts that any track theory is equivalent to a
strong one. An additive track theory is a track theory which moreover possesses
a lax zero object and
   finite
lax coproducts such that the natural map from the lax
coproduct to the lax product is a homotopy equivalence and the corresponding
homotopy category is an additive category. An additive track theory is called
very strong if it possesses a strong zero object and strong products which are also
strong coproducts. By the result of \cite{BJP}, any additive track theory is
equivalent to one which possesses strong products and lax coproducts or strong
coproducts and lax products. We show that in general it is impossible to get both
strong products and coproducts. However this is possible if certain obstructions
vanish. In particular this is possible if hom's of the corresponding homotopy
category are vector spaces over a field.
%Our interest with Shukla cohomology comes  with two points: firstly Shukla
%cohomology is  relation to the Mac Lane cohomology (see Chapter 13 of
%\cite{HC}) via spectral sequence \cite{PW}, \cite{P3}  and on the other hand
%Shukla cohomology in the dimension 3 classifies crossed bimodules (see Theorem
%\ref{shugaf}). We will show that very strongly additive track theories are
%intimately related with crossed bimodules, while  additive track theories are
%classified using third dimensional Mac Lane cohomology. Therefore the
%comparison between third dimensional Shukla and Mac Lane cohomology measures
%whether given an additive track theory is coming from a crossed bimodule or
%not.

 The contents of the sections below are as follows. In Section \ref{hochcoh}
 we recall basics on Hochschild cohomology theory and especially relationship
 between abelian extensions which are split over ground ring, and elements of
 the second Hochschild cohomology. In Section \ref{crossed} we
 introduce crossed bimodules and crossed extensions. We recall the relationship
  between crossed extensions which are split over the ground ring, and elements
of the third  Hochschild cohomology. This section also contains a new
interpretation of the classical obstruction theory in terms of crossed extensions
(see Theorem \ref{jvarediancinaag}). We also discuss a different
 generalization of the relationship between different sort of extensions and
 higher cohomology. In Section \ref{shuklacoh}
 we define Shukla cohomology as a kind of derived Hochschild cohomology on the
 category of chain algebras and we prove basic properties of the Shukla
 cohomology including relationship with crossed bimodules. In the original
 paper Shukla used an explicit cochain complex for the definition of Shukla
 cohomology. Unfortunately this complex is very complicated to work
 with. Quillen instead used closed model category structure on the category of
 simplicial algebras. We use the closed model category structure on the
 category of chain algebras, which is developed in the Appendix.
The Section \ref{shgam} is devoted to some computations of Shukla cohomology
when the ground ring is the ring of integers or $\Z/p^2\Z$; we also consider the
relationship between the Shukla cohomology over integers and over $\Z/p^2\Z$.
In this direction we prove the following result. Let $A$ be an algebra over
$\F_p$ and let $M$ be a bimodule over $A$, then the base change morphism
$$\Sh^i(A/\K,M)\to \Sh^i(A/\Z,M),\ \ \ \K=\Z/p^2\Z$$
is always an epimorphism. It is an  isomorphism in dimensions $0,1$ and $2$. We
also prove that in dimension three the kernel of this map is isomorphic to
$\H^0(A,M)$.  The Section \ref{relh} solves the problem of constructing a
canonical cochain complex for computing the Shukla cohomology in the important
case when the ground ring is an algebra over a field. Our cochain complex consists
of tensors, unlike the one proposed by Shukla. The Section \ref{hmlsh} recalls
basics of Mac Lane cohomology \cite{mac} and relationship with Shukla
cohomology. It is well known that these two theories are isomorphic up to dimension
two. It turns out that for algebras over fields they are also isomorphic in dimension three.
The section \ref{ota} continues the study of track theories which was started in \cite{BJP}.
In this section we show that the straightforward version of the strengthening result
for additive track theories is not true and we construct the corresponding obstruction.
This obstruction is defined using the exact sequence relating third Shukla and Mac Lane cohomology
and is a main application of the theory considered in the previous sections. The
Appendix contains the basic definitions on closed
model categories. It contains also a proof of the fact that chain algebras over any ground
ring form a closed model category. This fact is used in Section \ref{shuklacoh}.
At the end of the Appendix we introduce a closed model category structure on
the category of crossed bimodules over any ground ring.

In a forthcoming paper we introduce the notion of a {\it strongly additive track
theory} and we will prove that any additive track category is equivalent to a
strong one. The notion of strongly additive track theory is based on theory of
square rings \cite{square}.

 The second author is indebted to Mamuka Jibladze
for the idea to modify classical obstruction theory in terms of crossed
bimodules.

\section{ Preliminaries on Hochschild Cohomology}\label{hochcoh}
Here we recall the basic notion on  Hochschild cohomology theory
and refer to \cite{HC} and \cite{homology} for more details. In
this section $\K$ denotes a commutative ring with unit, which is
considered as a ground ring, except for the section \ref{relh}.
%Thus all modules are defined over
%$\K$ and we assume  that they are unitary, meaning that  the unite
%$1$ of $\K$ acts as the identity  operator and all algebras are associative
%$\K$-algebras with unites. A linear map means a $\K$-module homomorphism.
\subsection{Definition}
Let $R$ be a $\K$-algebra with unit and let $M$ be a bimodule
over $R$. Consider the module
$$
C^n(R,M):= \Hom(R^{\otimes n},M)
$$
(where $ \otimes =\otimes _K$ and $\Hom=\Hom_K$). The {\it
Hochschild coboundary} is the linear map $d: C^n(R,M)\to
C^{n+1}(R,M)$ given by the formula
\begin{multline*}
d(f)(r_1,...,r_{n+1})=r_1f(r_2,...,r_{n+1})+\\
+\sum_{i=1}^n(-1)^if(r_1,...,r_ir_{i+1},...,r_{n+1})+\\
(-1)^{n+1}f(r_1,...,r_n)r_{n+1}.
\end{multline*}
Here $f\in C^n(R,M)$ and $r_1,\cdots,r_{n+1}\in R$. By definition
the {\it $n$-th Hochschild cohomology group} of the algebra $R$
with coefficients in the $R$-bimodule $M$ is the $n$-th homology
group of the Hochschild cochain complex $C^*(R,M)$ and it is denoted by
$\H^*(R,M)$. Sometimes these groups are denoted by
$$\H^*(R/\K,M)$$
in order to make clear that the ground ring is $\K$. We are
especially interested in cases $\K=\Z,\F_p,\Z/p^2\Z$. It is clear
that for  such a $\K$ one has
%={\Z}$ or $\K={\Z}/p^k\Z$ one has
$$\H^i(\F_p/\K, \F_p)=0, \ i\geq 1.$$
In the following sections we consider two  modifications of Hochschild
cohomology, known as Shukla and Mac Lane cohomology. It should be
noted that in both theories the algebra $\F_p$ has nontrivial
cohomology over the ground ring $\K=\Z$ or $\K=\Z/p^2\Z$.

\subsection{$\K$-split exact sequences}
Let
$$\xymatrix{0\ar[r]&M_1\ar[r]^{\mu}&M\ar[r]^{\sigma}&M_2\ar[r]&0}$$
be an exact sequence of bimodules over $R$. It is called
\emph{$\K$-split} if there exists a $\K$-linear map $u:M_2\to M$ such
that $\sigma\circ u=\id_{M_2}$. This condition is equivalent to
the following one: there is a $\K$-linear map $v:M\to M_1$ such
that $v\circ \mu=\id_{M_1}$.

Let $f:M\to N$ be a morphism of bimodules over $R$. It is called
$\K$-split, if the following exact sequences
$$0\to \Ker(f)\to M\to {\sf Im}(f)\to 0,$$
and $$ 0\to {\sf Im}(f)\to N\to {\sf Coker}(f)\to 0$$ are
$\K$-split.

If
$\xymatrix{0\ar[r]&M_1\ar[r]^{\mu}&M\ar[r]^{\sigma}&M_2\ar[r]&0}$
is a $\K$-split exact sequence, then $$ 0\to C^*(R,M_1)\to
C^*(R,M)\to C^*(R,M_2)\to 0$$ is exact in the category of cochain
complexes and therefore yields the long cohomological exact
sequence:
$$\cdots \to \H^n(R,M_1)\to \H^n(R,M)\to \H^n(R,M_2)\to \H^{n+1}(R,M_1)\to
\cdots$$
\subsection{Induced bimodules} The category of
bimodules over $R$ is equivalent to the category of left modules
over the ring $R^e:=R\t R^{op}$, where $R^{op}$ is the opposite
ring of $R$, which is isomorphic to $R$ as a $\K$-module via the
map $r\mapsto r^{op}$, $R\to R^{op}$, while the multiplication
structure in $R^{op}$ is given by $r^{op}s^{op}=(sr)^{op}$. The
multiplication map $R\t R^{op}\to R$ is an algebra homomorphism.
We always consider $R$ as a bimodule over $R$ via this
homomorphism.

If $A$ and $B$ are left $R$-modules, then $\Hom(A,B)$ is a
bimodule over  $R$ by the following action
$$(rfs)(a)=rf(sa), \ \ r,s\in R, a\in M, f\in \Hom(A,B)$$
A bimodule is called {\it induced} if it is isomorphic to
$\Hom(R,A)$ for an $R$-module $A$. It is well-known
\cite{homology} that the Hochschild cohomology vanishes in
positive dimensions on induced bimodules.
 For a bimodule $M$ the map $$\mu:M\to \Hom(R,M)$$ given by
$\mu(m)(r)=mr$ is a homomorphism of bimodules, which is also
$\K$-split monomorphism, hence one has a $\K$-split short exact
sequence
$$0\to M\to \Hom(R,M)\to N\to 0$$
where $N={\sf Coker}(\mu)$, which yields the isomorphism
\begin{equation}\label{induced}
\H^{i+1}(R,M)\cong \H^{i}(R,N),   \ \ i>0
\end{equation}
This shows that there is
a natural isomorphism \cite{homology}
$$\H^*(R,M)\cong \Ext^*_{R^{e},\K}(R,M)$$
where subscript $\K$ indicates that $\Ext$-groups in question are defined in
the framework of relative homological algebra, where the proper class consists
of $\K$-split exact sequences. If $R$ is projective as a $\K$-module, then
one can take the usual $\Ext$-groups $\Ext^*_{R^{e}}(R,M)$  instead of the
relative $\Ext$-groups. In particular, the Hochschild cohomology vanishes in
positive dimensions on injective bimodules, provided $R$ is projective as a
$\K$-module.

\subsection{Hochschild cohomology in dimension $0$}
For $n=0$ one has
$$
\H^0(R,M)=\left\{ m\in M \mid rm=mr \ {\rm for \ any }\ r\in
R\right\} .
$$
In particular $\H^0(R,R)$ coincides with the center ${\sf Z}(R)$
of the algebra $R$.

\subsection{Hochschild cohomology in dimension $1$}
For $n=1$ a 1-cocycle is a linear map $D:R\to M$
satisfying the identity
$$
D(xy)=xD(y)+D(x)y, \ \ x,y\in R.
$$
Such a map is called a {\it derivation} from $R$ to $M$ and the
$\K$-module of derivations is denoted by $\der(R,M)$. A derivation $D:R\to M$
is a coboundary if it has the form $$\ad_m(r)=rm-mr$$ for some fixed
$m\in M$; $\ad_m$ is called an {\it inner derivation}. Therefore
$$
\H^1(R,M)=\der(R,M)/\lbrace {\rm Inner\ derivations} \rbrace .
$$
In particular one has the exact sequence
$$\xymatrix{0\ar[r]&\H^0(R,M)\ar[r]&M\ar[r]^-{\ad}\ &\ \der(R,M)\ar[r]&
\H^1(R,M)\ar[r]&0}$$
%The $A$-module $\H^1(R,R)$ is known as the
%$A$-module of outer derivations.

\subsection{Hochschild cohomology in dimension $2$}
It is clear that a $2$-cocycle of $C^*(R,M)$  is a linear map
$f:R\t R\to M$ satisfying
$$xf(y,z)-f(xy,z)+f(x,yz)-f(x,y)z=0, \ \ x,y,z\in R.$$
For any  linear map $g:R\to M$ the formula
$f(x,y)=xg(y)-f(xy)+f(x)y$ defines a cocycle, all such cocycles
are called coboundaries. We let ${\sf Z}^2(R,M)$ and ${\sf B}^2(R,M)$ be 
the collections  of all $2$-cocycles and coboundaries. Hence
$$\H^2(R,M)={\sf Z}^2(R,M)/{\sf B}^2(R,M).$$
We recall the relation of  $\H^2(R,M)$ to abelian extensions of
algebras.

 An {\it abelian extension} (sometimes called also {\it a singular  extension})
of associative
algebras is a short exact sequence
$$\xymatrix{
0\ar[r]& M\ar[r] &E\ar[r]^p& R\ar[r] &0} \leqno (E)
$$
where $R$ and $E$ are associative algebras with unit and $p:E\to
R$ is a homomorphism of algebras  with unit and $M^2=0$, in
other words the product in $E$ of any two elements from $M$ is
zero. For an elements $m\in M$ and $r\in R$ we put $mr:=me$ and
$rm=:em$. Here $e\in E$ is an element such that $p(e)=r$. This
definition does not depend on the choice of $e$. Therefore $M$
has a bimodule structure over $R$.

 An abelian  extension $(E)$ is called $\K$-split if there exists
 a linear map $s:R\to E$ such that $ps=\id_R$.

Assume we have a bimodule $M$ over an associative algebra $R$,
then we let ${\cal E}(R,M)$ be the category, whose objects are the
abelian extensions $(E)$ such that the induced $R$-bimodule
structure on $M$ coincides with the given one. The morphisms $(E)\to
(E')$ are commutative diagrams
$$
\xymatrix{0\ar[r] &M\ar[r]\ar[d]^{\id} &E\ar[r]\ar[d]^{\phi} &R\ar[r]\ar[d]^{\id}&0\\
         0\ar[r] &M\ar[r] &E'\ar[r] &R\ar[r]&0}
$$
where $\phi$ is a homomorphism of  algebras with unit. Moreover,
we let ${\cal E}_{\K}(R,M)$ be the category, whose objects are
$\K$-split singular extensions. It is clear that the categories
${\cal E}(R,M)$ and ${\cal E}_{\K}(R,M)$ are groupoids, in other
words all morphisms in ${\cal E}(R,M)$ and ${\cal E}_{\K}(R,M)$
are isomorphisms. We let ${\sf Extalg}(R,M)$ and  ${\sf
Extalg}_{\K}(R,M)$ be the classes of connected components of these
categories. Clearly ${\sf Extalg}_{\K}(R,M)\subset{\sf
Extalg}(R,M)$. According to \cite{homology} there is a natural
bijection
\begin{equation}\label{abgaf}
\H^2(R,M)\cong {\sf Extalg}_{\K}(R,M).
\end{equation}
We also recall that the map $\H^2(R,M)\to {\sf Extalg}_{\K}(R,M)$
is given as follows. Let $f:R\t R\to M$ be a $2$-cocycle. We let
$M\rtimes _f R$ be an associative $\K$-algebra which  is $M\oplus
R$ as a $\K$-module, while the algebra structure is given by
$$(m,r)(n,s)=(ms+rn+f(r,s),rs).$$
Then
$$\xymatrix{0\ar[r]&M\ar[r]^-i \ & \ M\rtimes _fR\ \ar[r]^-p\ &R\ar[r]&0}$$
is an object of ${\cal E}_{\K}(R,M)$. Here $i(m)=(m,0)$ and
$p(m,r)=r$.

\subsection{Cohomology of tensor algebras.}  \cite{HC}, \cite{homology}. Let
$V$ be a $\K$-module. For the tensor algebra $R=T^*(V)$ one has
$H^i(R,-)=0$ for all $i\geq 2$. An algebra is called {\it free} if
it is isomorphic to $T(V)$, where $V$ is a free $\K$-module.

\subsection{Cup-product in Hochschild cohomology} For any
associative algebra $R$ the cohomology $\H^*(R,R)$ is a graded
commutative algebra under the cup-product, which is defined by
$$(f\cup g)(r_1,\cdots
,r_{n+m}):=f(r_1,\cdots,r_n)g(r_{n+1},\cdots,r_{n+m})$$, for $f\in
C^n(R,R)$ and $g\in C^m(R,R)$ (see \cite{gerstenhaber}). This
product corresponds to the Yoneda product under the isomorphism
$\H^*(R,R)\cong \Ext^*_{R^{e},\K}(R,R)$.

%%%%%%%%%%%%%%%%%%%%%%%%%%%%%%%%%%%%%%%%%%%%%%%%%%%%%%%%%%%%%%%%%%%%%%
%%%%%%%%%%%%%%%%%%%%%%%%%%%%%%%%%%%%%%%%%%%%%%%%%%%%%%%%%%%%%%%%%%%%%
%%%%%%%%%%%%%%%%%%%%%%%%%%%%%%%%%%%%%%%%%%%%%%%%%%%%%%%%%%%%%%%%%%%%
\section{Crossed bimodules and Hochschild cohomology }\label{crossed}
\subsection{Crossed bimodules}

Let us recall that a chain algebra over $\K$ is a graded algebra
$C_*=\bigoplus_{n\geq 0} C_n$ equipped with a boundary map
$\partial: C_*\to C_*$ of degree $-1$ satisfying the Leibniz
identity
$$\partial (xy)=\partial (x)+(-1)^{|x|}x \partial(y).$$
\begin{De} A {\it  crossed bimodule} is a chain algebra which is trivial  in
dimensions $\geq 2$.
\end{De}

Thus a crossed bimodule consists of an algebra $C_0$ and a bimodule $C_1$ over $C_0$ together with
a homomorphism of bimodules
$$C_1 \buildrel \partial \over \to C_0$$
such that
$$\partial (c)c'=c\partial (c'), \ c,c'\in C_1,$$
Indeed, since $C_2=0$ the last condition is equivalent to the  Leibniz identity
$0=\partial(cc')=\partial(c)c'-c\partial(c')$.

It follows that the product  defined by $$c*c':=\partial (c)c'$$
where $c,c'\in C_1$ gives an associative non-unital $\K$-algebra
structure on $C_1$ and $\partial:C_1\to C_0$ is a homomorphism of
non-unital $\K$-algebras. The equivalent but less economic
definition goes back at least to Dedecker and Lue \cite{delu}. The
notion of crossed bimodules is an associative algebra analogue of
crossed modules introduced by Whitehead \cite{hanry} in the group
theory framework, which plays a major role in homotopy theory of
spaces with nontrivial fundamental groups \cite{B1}, \cite{jll}.

We let ${\sf Xmod}$ and ${\sf Xmod}_R$
be the category of crossed bimodules and crossed $R$-bimodules
respectively.

 We have also a  category
 ${\sf Bim/Alg}$, whose objects are triples $(C_0,C_1,\d)$,
where $C_0$ is an associative algebra, $C_1$ is a bimodule over
$C_0$ and $\d:C_1\to C_0$ is a homomorphism of bimodules over
$C_0$. It is clear that ${\sf Xmod}$ is a  full subcategory of
${\sf Bim/Alg}$ and the inclusion  ${\sf Xmod}\subset {\sf
Bim/Alg}$ has a left adjoint functor, which assigns $\d:{\wt
C_1}\to C_0$ to $C_1\to C_0$. Here ${\wt C_1}$ is the quotient of
$C_1$ under the equivalence relation $x\d(y)-\d(x)y\sim 0$,
$x,y\in C_1$.

We let ${\sf Mod/Alg}$ be the category whose objects are triples $(V,C,\d)$,
where $C$ is an associative algebra, $V$ is a $\K$-module and $\d:V \to C$ is a
linear map. One has the forgetful functor ${\sf Bim/Alg}\to {\sf Mod/Alg}$,
which has a left adjoint functor sending $(V, C,\d)$ to the triple
$(M,d,C)$, where $M=C\t V\t C$ and $d$ is the unique homomorphism of bimodules
which extends $\d$. As a consequence we see that the forgetful functor ${\sf
Xmod}\to {\sf Mod/Alg}$ also has a left adjoint.
%, whose value
%on$(V,\d,C)$ is called {\it free crossed $C$-bimodule generated by $d:V\to C$.}
Of special interest is the case when $C$ is a free associative algebra and
$V$ is a free $\K$-module on $X \subset V$. In this case the corresponding
crossed bimodule is called  {\it free crossed bimodule}.

\subsection{Hochschild cohomology in the dimension $3$ and crossed extensions}
Here we recall the relation between Hochschild cohomology and crossed
bimodules (see Exercise  E.1.5.1 of \cite{HC}  or
\cite{baumin}).

\smallskip
Let $\partial: C_1\to C_0$ be a crossed bimodule. We put $M=\ker
(\partial)$ and $R=\Cok (\partial)$. Then the image of  $\partial$
is an ideal of $C_0$. We have also $MC_1=0=C_1M$ and $M$ has  a
well-defined bimodule structure over $R$.

Let $R$ be an associative algebra with unit and let $M$ be a
bimodule over $R$. A {\it  crossed extension} of $R$ by $M$ is an
exact sequence
$$0\to M\to C_1 \buildrel \partial \over \to C_0\to R\to 0$$
where $\partial:C_1\to C_0$ is a crossed bimodule,  such that
$C_0\to R$ is a homomorphism of algebras with unit and an
$R$-bimodule structure on $M$ induced from the crossed bimodule
structure coincides with the prescribed one.

A crossed extension of $R$ by $M$ is {\it $\K$-split}, if all
arrows in the exact sequence
$$0\to M\to C_1 \buildrel \partial \over \to C_0\to R\to 0$$
are  $\K$-split.

For fixed $R$ and $M$ one can consider the category ${\bf
Crossext}(R,M)$ whose objects are crossed extensions of $R$ by
$M$. Morphisms are maps between crossed modules which induce the
identity on $M$ and $R$.

\begin{Le}\label{crospullback}
Assume $(\d)$ is a crossed extension of $R$ by $M$ and  a
homomorphism $f:P_0\to C_0$ of unital $\K$-algebras is given. Let
$P_1$ be the pull-back of the diagram
$$\xymatrix{&P_0\ar[d]\\
C_1\ar[r]&C_0.}$$ Then there exists a unique crossed module
structure on $P_1\to P_0$ such that the diagram
$$\xymatrix{0\ar[r]&M\ar[d]^{\id}\ar[r]&P_1\ar[r]\ar[d]&P_0\ar[d]^f
\ar[r]&R\ar[r]\ar[d]^{\id}&0\\
0\ar[r] & M\ar[r] &C_1\ar[r]&C_0\ar[r] &R\ar[r] &0}$$ defines a
morphism of  crossed extensions.
\end{Le}

\begin{Co}\label{crosfiltri}
In each connected component of
${\bf Crossext}(R,M)$ there is a crossed extension $$0\to M\to P_1
\to P_0\to R\to 0 \leqno(P)$$ with free algebra $P_0$ and for any
other object $(\d)$ in this connected component
there is a morphism $(P)\to (\d)$. Thus  $(\d)$ and $(\d')$ are in
the same component of ${\bf Crossext}(R,M)$ iff there exists a
diagram of the form $(\d)\leftarrow (P)\to (\d')$.
\end{Co}

We let ${\bf Crossext}_{\K}(R,M)$ be the subcategory of $\K$-split
crossed extensions. Morphisms are such morphisms from ${\bf
Crossext}(R,M)$ that all maps involved are $\K$-split. Let ${\bf
Cros}(R,M)$ and ${\bf Cros}_{\K}(R,M)$ be the set of components of
the category of crossed extensions and the category of $\K$-split
crossed extensions respectively. Then there is a canonical
bijection:
\begin{equation}\label{bauh3}
\H^3(R,M)\cong {\bf Cros}_{\K}(R,M)
\end{equation}
(see for example Exercise E.1.5.1 of \cite{HC} or \cite{baumin}).
A similar isomorphism for group cohomology was proved by Loday
\cite{jllst}, see also \cite{MW}.
 We recall only how to associate  a 3-cocycle to
a $\K$-split crossed extension:
$$0\to M\to C_1 \buildrel \partial \over \to
C_0\buildrel \pi \over \to R\to 0$$
of $R$ by $M$.
We put $V:={\sf Im}(\partial)$ and consider $\K$-linear sections $p:R\to C_0$
and $q:V\to C_1$
of $\pi: C_0\to R$ and $\partial:C_1\to V$ respectively. Now we define
$m:R\t R\to V$
by $m(r,s):=q(p(r)p(s)-p(rs))$. Finally we define
$f:R\t R\t R\to M$ by
$$f(r,s,t):=p(r)m(s,t)-m(rs,t)+m(r,st)-m(r,s)p(t).$$
Then $(f,g,h)\in {\sf Z}^3(R/\K,M)$ and the corresponding class in $\H^3(R/\K,M)$
depends only on the
connected component of a given crossed extension  and in this way one gets the expected isomorphism (see \cite{baumin}).

\subsection{Obstruction theory}
 Now we explain a variant of the classical obstruction theory in
terms of crossed extensions (compare with Sections IV.8 and IV.9 of
\cite{homology}). Let
$$0\to M\to C_1 \buildrel \partial \over \to C_0\to R\to 0  $$
be a crossed extension of $R$ by $M$.
{\it A $\d$-extension of $C_1$ by $R$}
is a commutative diagram with exact rows
$$\xymatrix{&0\ar[r]&C_1\ar[d]^{\id}
\ar[r]^{\mu}&S
\ar[d]^{\xi}\ar[r]^{\si}
&R\ar[d]^{\id}\ar[r]&0\\
0\ar[r]&M\ar[r]&C_1\ar[r]^\d&C_0\ar[r]&R\ar[r]&0}$$
where $S$ is a unital $\K$-algebra and $\si$ is a
 homomorphism of unital $\K$-algebras. Furthermore
 the equalities
 $\mu(x)s=\mu(x\xi(s)))$
 and $s\mu(x)=\mu(\xi(s)x)$
 hold,
where $x\in C_1$, $s\in S$. It follows then that product in $C_1$
induced from $S$ coincides with the $*$-product:
$x*y=\d(x)y=x\d(y)$. Moreover one has the exact sequence
$$\xymatrix{0\ar[r]& M\ar[r]^{\mu}
& S\ar[r]^{\xi}& C_0\ar[r]& 0.}$$
It is clear that $\d$-extensions of $C_1$ by $R$ form a groupoid, whose set of components will be denoted by
 $_\d{\Ext}(R,C_1)$.

Now we assume that $\d$ is  a $\K$-split crossed
extension. A $\d$-extension of $C_1$ by $R$ is called $\K$-split if $\xi$
 is a $\K$-split epimorphism. Of course in this case $\si$ is $\K$-split as well.
 We let $_\d\Ext_{\K}(R,C_1)$ be the subset
 of $_\d{\Ext}(R,C_1)$ consisting of $\K$-split $\d$-extensions.
\begin{The}\label{jvarediancinaag} The class of a $\K$-split crossed extension
$$0\to M\to C_1 \buildrel \partial \over \to C_0\to R\to 0 \leqno(\d)$$
is zero in $\H^3(R,M)$  iff $_\d\Ext_\K(R,C_1)$ is nonempty. If this
is the case then the group $\H^2(R,M)$ acts transitively and effectively
on $_\d\Ext_\K(R,C_1)$.
\end{The}
{\it Proof}. For a crossed extension $\d$ one considers sections
$p:R\to C_0$ and $q:V\to C_1$, $V= {\sf Im}(\d)$
 as above. We may and we will assume that $p(1)=1$.
Then the class
of $(\d)$
 in $\H^3$ is given by the cocycle
$f(r,s,t):=p(r)m(s,t)-m(rs,t)+m(r,st)-m(r,s)p(t)$
where $m(r,s)=q(p(r)p(s)-p(rs))$.
Given a $\d$-extension of $C_1$ by $R$:
$$\xymatrix{&0\ar[r]&C_1\ar[d]^{\id}
\ar[r]^{\mu}&S
\ar[d]^{\xi}\ar[r]^{\si}
&R\ar[d]^{\id}\ar[r]&0\\
0\ar[r]&M\ar[r]&C_1\ar[r]^\d&C_0\ar[r]&R\ar[r]&0}$$ choose a
$\K$-linear section $v:C_0\to S$ such that $v(1)=1$.
  One puts
$u=vp:R\to S$. Then $\si u= {\sf Id}_R$. One defines $n:R\t R\to
C_1$ by $\mu (n(r,s))=u(r)u(s)-u(rs)$. We claim  that
\begin{equation}\label{cinakocikli}
p(r)n(s,t)-n(rs,t)+n(r,st)-n(r,s)p(t)=0
\end{equation}
Indeed, $$p(r)n(s,t)=u(r)n(s,t)=u(r)u(s)u(t)-u(r)u(st).$$
Similarly $n(r,s)p(t)=u(r)u(s)u(t)-u(rs)u(t)$. Thus
$$p(r)n(s,t)-n(rs,t)+n(r,st)-n(r,s)p(t)=u(r)u(s)u(t)-u(r)u(st)-u(rs)u(t)$$
$$+u(rst)+u(r)u(st)-u(rst)-u(r)u(s)u(t)+u(rs)u(t)=0.$$
Since $m(r,s)=q\d n(r,s)$, it follows
that $g(r,s)=m(r,s)-n(r,s)$ lies in $M$. Thus we obtain
a well-defined linear map $g:R\t R\to M$. Then it follows from
the equation (\ref{cinakocikli})
that
$$f(r,s,t)=rg(s,t)-g(rs,t)+g(r,st)-g(r,s)t,$$
which shows that the class of $\d$
 in $\H^3$ is  zero. Given
any normalized $2$-cocycle $h:R\t R\to M$, one can define a new
$\d$-extension $S_h$ of $R$ by $C_1$ by putting $S_h=C_1\os R$
with the  following  multiplication:
$$(x,r)(y,s)=(x*y+p(r)y+xp(s)+n(r,s)+h(x,y),xy).$$
This construction yields a transitive and effective action of
$\H^2\!(R,\!M)$ on $_\d\!\Ext_\K\!(R,\!C_1)$.

Conversely, assume that the class of $0\to M\to C_1 \buildrel
\partial \over \to C_0\to R\to 0$ is zero in $\H^3(R,M)$. Thus
there is a linear map $g:R\t R\to M$ such that
$f(r,s,t)=rg(s,t)-g(rs,t)+g(r,st)-g(r,s)t.$ One can define $n: R\t
R\to C_1$ by $n(r.s)=m(r,s)-g(r,s)$. Then
$p(r)n(s,t)-n(rs,t)+n(r,st)-n(r,s)p(t)=0$ and therefore $S=R\os
C_1$ with the product $(x,r)(y,s)=(x*y+p(r)y+xp(s)+n(x,y),xy)$
defines a $\d$-extension.

\subsection{ Abelian and crossed $n$-fold extensions}\label{baumin-n}
{\it An abelian twofold extension} of an algebra $R$ by an
$R$-$R$-bimodole $M$ is an exact sequence
$$\xymatrix{0\ar[r] &M\ar[r]^{\aa}& N\ar[r]^{\mu} &S\ar[r]^{\pi}& R\ar[r] &0}$$
where $N$ is a bimodule over $R$ and $\aa$ is a bimodule
homomorphism. Moreover, $S$ is an associative algebra with unit
and $\pi$ is a homomorphism of algebras with unit, such that
$\Ker(\pi)$ is a square zero ideal of $S$. Furthermore, for any
$s\in S$ and $n\in N$ one has
$$\mu(n\pi(s))=\mu(n)s, \ \ \mu(\pi(s)n)=s\mu(n).$$

 We let ${\cal E}^2(R,M)$ be  the category of
abelian twofold extensions of $R$ by $M$, whose connected
components are denoted by ${\bf Extalg}^2(R,M)$. As usual we have
also a $\K$-split variant ${\cal E}^2_\K(R,M)$ of the category
${\cal E}^2(R,M)$: Objects of ${\cal E}^2_\K(R,M)$ are $\K$-split
twofold abelian extensions (i.e. $\aa$, $\mu$ and $\pi$ are
$\K$-splits), and  the morphisms in ${\cal E}^2_\K(R,M)$ are
$\K$-splits, accordingly we let ${\bf Extalg}^2_\K(R,M)$ be the
connected components of ${\cal E}^2_\K(R,M)$.

Let us note that for any abelian twofold extension
$$\xymatrix{0\ar[r] &M\ar[r]^{\aa}& N\ar[r]^{\mu} &S\ar[r]^{\pi}& R\ar[r] &0}$$
the morphism $\mu:N\to S$ is a crossed bimodule, where the action
of $S$ on $N$ is given via $\pi$. It is clear that the induced
$*$-product on $N$ is trivial. Thus one obtains the functor ${\cal
E}^2(R,M)\to {\bf Crossext}(R,M)$, which takes the subcategory
${\cal E}^2_\K(R,M)$ to the category ${\bf Crossext}_\K(R,M)$.
\begin{Le}\label{2abgaf}
The natural map ${\bf Extalg}^2_\K(R,M)\to {\bf Cros}_\K(R,M)$
is a bijection and therefore
$${\bf Extalg}^2_\K(R,M)\cong\H^3(R,M)$$
\end{Le}

{\it Proof}. We  just construct the inverse map
$$\xi:\H^3(R,M)\to {\bf Extalg}^2_\K(R,M).$$
Consider the $\K$-split short exact sequence
$$0\to M\to \Hom(R,M)\to N\to 0$$
and the corresponding  isomorphisms (\ref{induced}), (\ref{abgaf})
$$\H^3(R,M)\cong \H^2(R,N)\cong {\bf Extalg}(R,N).$$ Take now an
element $a\in \H^3(R,M)$. It corresponds under these isomorphisms
to an abelian extension $0\to N\to S\to R\to 0$. By gluing it with
$0\to M\to \Hom(R,M)\to N\to 0$ one obtains an abelian twofold
extension
$$0\to M\to \Hom(R,M)\to S\to R\to 0$$
In this way one obtains the expected map $\xi$.

It is clear now  how to introduce the notion of abelian $n$-fold
extension for all $n\geq 2$ and get the same sort of isomorphism
in higher dimensions.

Following the earlier work of Huebschmann \cite{hueb},
recently Baues and Minian \cite{baumin} obtained another
interpretation of Hochschild cohomology in dimensions $\geq 4$.
They introduced the notion of crossed $n$-fold extension and
proved that $n$-fold extensions classify $(n+1)$-dimensional
Hochschild cohomology for all $n\geq 2$. For $n=2$ this is an
isomorphism (\ref{bauh3}). Here we give  a sketch how to deduce
the case $n>2$ from the case $n=2$ and  from the classical results
of Yoneda \cite{yoneda}. This argument gives also a new proof of
Lemma \ref{2abgaf}.

Let $T$ be an additive functor from the category of bimodules over
$R$ to the category of $\K$-modules. Objects of the category
${\cal E}^n(T)$ are pairs $(E,x)$, where
$$0\to M\to E_1\to \cdots \to E_n\to 0$$
is a $n$-fold extension of $E_n$ by $M$ in the category of
$R$-$R$-bimodules and $x\in T(E_n)$. Morphisms in ${\cal E}^n(T)$
are defined in an obvious way. Let ${\bf E}^n(T)$ be the set of
components of the category ${\cal E}^n(T)$. A result of Yoneda
asserts that one has a natural isomorphism:
$${\bf E}^n(T)\cong S^n T(M)$$
where $S^nT$ is the $n$-th satellite of $T$ \cite{CE}.

Comparing with the definition of abelian twofold extension we see
that
$${\bf Exalg}^2(R,M)\cong {\bf E}^1(T)$$ where $T={\bf
Extalg}(R,-)$. To show how to deduce Lemma \ref{2abgaf} from the
Yoneda isomorphism, we consider the case when $\K$ is a field.
Since $T\cong \H^2(R,-)$ the result of Yoneda yields
$${\bf Extalg}^2(R,M)\cong S^1T(M)\cong \H^3(R,M).$$ This argument works also for  general
$\K$: we have to use a straightforward generalization of the
Yoneda isomorphism in the framework of relative homological
algebra.

Let $n\geq 2$. A {\it crossed $n$-fold extension of $R$ by $M$}
\cite{baumin} is an exact sequence
$$\xymatrix{0\ar[r]&M\ar[r]^-{f}&M_{n-1}\ar[r]^{\d_{n-1}}&\cdots \ar[r]^{\d_2}&
M_1\ar[r]^{\d_1}& A\ar[r]^{\pi}& R\ar[r]&0 }
$$
of $\K$-modules with the following properties:

i) $(M_1,R,\d_1)$ is a crossed bimodule with cokernel $R$;

ii) $M_i$ is a bimodule over $R$ for $1<i\leq n-1$ and $\d_i$ and
$f$ are maps of bimodules  over $R$. Note that $\Ker(\d_1)$ is
naturally a bimodule over $R$ and therefore it makes sense to
require $\d_2$ to be a map of bimodules over $R$. We let ${\bf
Cros}^n(R,M)$ denote the set of connected components of the category of crossed
$n$-fold extensions of $R$ by $M$. Observe that
$${\bf Cros}^n(R,M)={\bf E}^{n-2}(T)$$
where $T={\bf Cros}(R,-)$.

Now, as in \cite{baumin} for simplicity we assume that $\K$ is a
field. Theorem 4.3 of \cite{baumin} claims that there is a natural
isomorphism
\begin{equation}\label{baun}
{\bf Cros}^n(R,M)\cong \H^{n+1}(R,M)
\end{equation}
 For $n=2$ this is the
isomorphism (\ref{bauh3}) and for $n>2$ it is an immediate
corollary of Yoneda's isomorphism:

\begin{multline*}
{\bf
E}^{n-2}(T)=S^{n-2}T=S^{n-2}\H^3(R,-)(M)=\\
S^{n-2}S^3\H^0(R,-)(M)\cong S^{n+1}\H^0(R,-)(M)=\H^{n+1}(R,M)
\end{multline*} Here we used the isomorphism $T\cong H^3(R,-)$ and
the classical fact that $\H^n(R,M)=\Ext_{R\t
R^{op}}^n(R,M)=S^n\H^0(R,-)(M)$ see \cite{CE}. For general $\K$
one needs to work in the framework of relative homological algebra
\cite{homology}. The corresponding class of proper exact sequences
consists of $\K$-split exact sequences. Then the corresponding
results hold for arbitrary $\K$.

 As we can see the
results in this section strongly depend on the vanishing of
Hochschild cohomology on (relative) injective modules.

\section{Shukla Cohomology}\label{shuklacoh}
As we already saw the Hochschild cohomology in dimensions two and three
classifies $\K$-split abelian and crossed extensions respectively. However,
there is a variant of Hochschild cohomology due to Shukla in the early 60-s
which classifies all abelian and crossed extensions. We will present these
results. Our approach to Shukla cohomology is based on chain algebras and
especially on the possibility of extension of Hochschild cohomology to chain
algebras. Actually there are two ways for such extension. First is a very naive:
one replaces $\otimes $ and $\Hom$ in the definition of Hochschild cohomology
by the tensor product and Hom of complexes to arrive at a cosimplicial cochain
complex and then one takes the homology of the total complex. However, this
definition does not respect weak equivalences of chain algebras. The second
definition (called derived Hochschild cohomology) is a kind of Quillen's
derivative of the Hochschild cohomology and uses
 the closed model category structure on the category of chain complexes
introduced in the Appendix. Since the category of algebras is the
 full subcategory of the category of chain algebras, the derived Hochschild
cohomology restricts to a cohomology theory of algebras, which is by definition the Shukla cohomology.

\subsection{Hochschild cohomology for chain algebras} In this section we give
a naive definition of the
Hochschild cohomology for chain algebras.

Let us recall that a chain algebra is a graded algebra
$R_*=\bigoplus_{n\geq 0} R_n$ equipped with a differential $d:R_n\to
R_{n-1}$ satisfying the Leibniz identity:
$$d(xy)=d(x)y+(-1)^{n}xd(y),   \ \ x\in A_n, y\in A_m.$$
Let ${\bf DGA}$ be the category of chain algebras. A morphism of
chain algebras is a weak equivalence if it induces isomorphism in homology.

An  $R_*$-bimodule  is a chain complex $M_*$ of $\K$-modules,
equipped with actions from both sides: $R_*\t M_*\to M_*$ and $M_*\t
R_*\to M_*$, satisfying usual axioms. However, for our purposes we
restrict ourselves to the case when $M$ is concentrated in
degree zero. In this case $R_*$-bimodule means simply a bimodule
over $H_0(R_*)$. In particular $xm=0=mx$ as soon as $m\in M$ and
$|x| \geq 1$. For a chain algebra $R_*$ and a
$H_0(R_*)$-bimodule $M$ we denote by $C_*(R_*,M)$ the total
complex of the following cosimplicial cochain complex. The $n$-th
component of this cosimplicial object is the cochain
complex
$$
C^n(R_*,M):= \Hom(R^{\otimes n}_*,M).
$$
Here $ \otimes$ denotes the tensor product of chain complexes. The
coface operations are given via Hochschild coboundary formula:
$$
d^0(f)(r_1,...,r_{n+1})=(-1)^{nk}r_1f(r_2,...,r_{n+1}), \  \ f:R^{\otimes n}_*
\to M, \  \ r_1\in R_k$$
(actually this expression is zero provided $k>0$)
$$
 d^i(f)(r_1,...,r_{n+1})= f(r_1,...,r_ir_{i+1},...,r_{n+1}), \  \ 0<i<n+1 $$
$$d^{n+1}(f)(r_1,...,r_{n+1})=f(r_1,...,r_n)r_{n+1}.$$
The homology of $C_*(R_*,M)$ is denoted by $\H^*(R_*,M)$.

The spectral sequences of a bicomplex in our situation have the
following form:
$$E_{pq}^1=H^q(\Hom(R_*^{\t p},M))\then \H^{p+q}(R_*,M)$$
$$F_{pq}^1=\H^q(| R_*|,M)\then\H^{p+q}(R_*,M)$$
Here $| R_*|$ denotes the underlying graded algebra of the chain algebra $(R_*,\d)$.

\begin{Le}\label{hhisinvariantoba} Let $f:R_*\to S_*$ be a weak equivalence of chain algebras
and let $M$ be a bimodule over $H_0(S)$. Then the induced homomorphism $$\H^*(S_*,M)\to \H^*(R_*,M)$$ 
s an isomorphism provided $R_*$ and $S_*$ are projective $\K$-modules.
%ap of chain algebras,
%hich induces isomorphism in homology in dimensions $<$.
%nd an epimorphism in dimension $m$
%nd let
%M$ be a bimodule over $H_0(S)$. Then the induced homomorphism
%\H^iS_*,M)\to \H^iR_*,M)$ is an isomorphism provided $R_*$ and $S_*$ are
%rojective $\K$-module and $i\leq m$.
\end{Le}

%\it Proof}.
%t follows from the K\"unneth and the universal coefficient
%heorems that the induced map
%\Hom(R_*^{\t n},M)
%to \Hom(S_*^{\t n},M)$ is also an
%somorphism on the homological level
%n dimensions $\leq n$. On the other hand on eof the
%he spectral sequence of bicomplex in our sitiation has th
{\it Proof}.
It is well known that any weak equivalence between  degreewise
projective bounded below chain complexes is a homotopy equivalence. Thus
$f$ is a homotopy equivalence in the category
of chain complexes of $\K$-modules. Therefore the induced map
$R_*^{n \t}\to S_*^{n \t}$ is also a homotopy equivalence. It follows that
the induced map of cosimplicial cochain complexes is a homotopy equivalence
in each degree and therefore it induces  a weak equivalence on the total
complex level thanks to the spectral sequence argument associated to the
 the bicomplex.
\subsection{The complex $\Der(| R_*|,M)$}\label{dercompleksi}
Let $R_*$ be a chain algebra and $M$ be an $H_0(R_*)$-bimodule. We
can take the derivations $| R_*|\to M$ from the underlying
graded algebra to $M$. Since $| R_*|$ is graded, the space
of derivations $\Der(| R_*|,M)$ is also graded. Since $M$ is
concentrated only in dimension zero, we see that the $0$-th
component is the space of all  derivations  $R_0\to M$, while in
dimensions $n>0$ we get the space of linear maps $f:R_n\to M$
satisfying the conditions
$$f(xy)=xf(y), \ \ f(yx)=f(x)y, \ \ x\in R_0, y\in R_n, n>0$$
$$f(uv)=0, \ \ u\in R_i, v \in R_j, i+j=n, i>0, j>0.$$
The boundary  map $\partial: R_n\to R_{n-1}$ in $R_*$ yields a
cochain complex structure on $\Der(| R_*|,M)$.
In what follows $\Der(| R_*|,
,M)$ is always considered with this cochain
complex structure.

\Lem{Let $R_*$ be a quasi-free algebra, meaning that the underlying algebra
structure is free, and let  $M$ be an $H_0(A_*)$-bimodule.
Then the Hochschild cohomology  $\H^n(A_*,M)$ is isomorphic to the
$(n-1)$-st homology of the cochain complex $\Der(| R_*|,M)$
provided $n>0$.
}

{\it Proof}. This is a direct consequence of the spectral sequence
related to the bicomplex $C^*(R_*,M)$ together with the fact that
the Hochschild cohomology of a free algebra is zero in dimensions
$>1$.

We also recall the  K\"unneth formula for Hochschild cohomology

\begin{Le}\cite{homology}\label{kiuneti}
 Let $\K$ be a field and let $R_*$ and $S_*$ be chain algebras. Assume that for each
 $n$, $R_n$ and $S_n$ are finite dimensional vector spaces. Then for any $R_*$-bimodule $M$ and $S_*$-bimodule
$N$ one has the following isomorphism
$$\H^n(R\t S, M\t N)\cong \bigoplus_{i+j=n}\H^i(R,M)\t \H^j(S,N)$$
\end{Le}
\subsection{Derived Hochschild cohomology and Shukla cohomology}
In this section we use the closed model category structure  on chain algebras
described in the Appendix. Let us recall that weak equivalences in this
model category are usual ones and   a morphism of chain algebras is a fibration
if it is surjective in all positive dimensions. We also need the fact that any
cofibrant object is a retract of a quasi-free algebra. It follows from Lemma
\ref{hhisinvariantoba} that for any weak equivalence $f:R_*\to S_*$ of
cofibrant chain algebras and any $H_0(S)$-bimodule $M$ the induced homomorphism
$\H^*(S_*,M)\to \H^*(R_*,M)$ is an isomorphism. We can use this fact to define
the derived Hochschild cohomology as follows. Let $R_*$ be a chain algebra.
Thanks to the properties of closed model categories there exists a chain algebra
morphism $f:R_*^c\to R_*$ which is a weak equivalence and $R_*^c$ is a
cofibrant. For any $R$-bimodule $M$ the groups $\H^*(R^c_*,M)$ do not depend
on the cofibrant replacement and they are called {\it the derived Hochschild
cohomology of $R_*$ with coefficients in $M$} and are denoted by ${\bf
H}^*(R_*,M)$. Thus
$${\bf H}^*(R_*,M):=\H^*(R^c_*,M)$$
This definition has expected functorial properties: for any morphism $f:R_*\to S_*$
of chain algebras and any $H_0(S)$-bimodule $M$ there is a well-defined
homomorphism ${\bf H}^*(S_*,M)\to {\bf H}^*(R_*,M)$ which depends only on the
homotopy class of $f$. Moreover it is an isomorphism provided $f$ is a weak
equivalence. One has also a natural homomorphism $\H^*(R_*,M)\to {\bf
H}^*(R_*,M)$ which is induced by the chain algebra homomorphism $R_*^c\to R_*$.
The following fact is a direct consequence of Lemma \ref{hhisinvariantoba}.
\begin{Le} If $R_*$ is projective as a $\K$-module then $\H^*(R_*,M)\to
{\bf H}^*(R_*,M)$ is an isomorphism.
\end{Le}

Since the category of algebras is a full subcategory of the category of chain
algebras we can consider the restriction of the derived Hochschild cohomology
${\bf H}^*$ on the category of algebras. The resulting theory is called {\it
the Shukla cohomology.} Thus for any algebra $R$ and any $R$-bimodule $M$ the
Shukla cohomology of an algebra $R$ with coefficients in $M$ is defined by
$$\Sh^*(R,M):={\bf H}^*(R,M)\cong H^*(\Der(R^c_*,M))$$ where $R^c_*\to R$ is a weak
equivalence from a quasi-free chain algebra $R^c_*$. The natural transformation
$$\H^n(R,M)\to \Sh^n(R,M), \ \ n\geq 0$$
is an isomorphism in dimensions $n=0,1$ and it is an isomorphism
in all dimensions provided $R$ is projective as a $\K$-module. For
example we have $\Sh^i(A,-)=0$ provided $A$ is a free algebra and $i\geq 2$.

The cup-product in Hochschild cohomology yields a (commutative graded)
algebra structure on $\Sh^*(A,A)$.

\subsection{Shukla cohomology and extensions}\label{gafshu}

The following properties of Shukla cohomology are of special interests. They
are non-$\K$-split analogues of the isomorphisms
 (\ref{abgaf}) and (\ref{bauh3}).

\begin{The}\label{shugaf} Let $A$ be an associative algebra and let $M$ be an $A$-bimodule.
Then there are natural isomorphisms
$$\Sh^2(A,M)\cong {\bf Extalg}(A,M)$$
$$\Sh^3(A,M)\cong {\bf Cros}(A,M).$$
\end{The}
The first isomorphism is well known (see Theorem 4 of \cite{sh}).
However we give an independent proof.

{\it Proof}. i) Let
 $$
0\to M\to E\to A\to 0 \leqno (E)
$$
be a singular extension of algebras. Define the chain algebra $E_*$ as
follows:
$$E_0=E, \ \ E_1=M, \ \ E_n=0, \ \ n\geq 2$$
The only nontrivial boundary map is induced by the inclusion $M\to
E$. Then one has a map of chain algebras $E_*\to A$ which is an
acyclic fibration. Let $A_*\to A$ be a weak equivalence with
quasi-free $A_*$. Since $A_*$ is cofibrant  there exists a
lifting $f_*:A_*\to E_*$. We consider now the first component
$f_1: A_*\to E_1=M$ of $f_*$. Since $f_*$ is a homomorphism of
algebras it follows that $f_1\in {\sf Der}(A_*,M)$ is a 1-cocycle
of $ {\sf Der}(A_*,M)$
  and therefore it gives  a class
$e(E)\in \Sh^2(A,M)$. If $g_*:A_*\to E_*$ is another lifting, then the values of
$h=f_0-g_0:A_0\to E$ lie in $M$. Thus $h\in {\sf Der}(A_0,M)$ and $f_1-g_1=
\partial ^*(h)$, which shows that the class $e(E)$ depends only
on the isomorphism class of $(E)$. Conversely, if
$f\in {\sf Der}(A_*,M)$ is a 1-cocycle, then one can form
an abelian extension
according to the following diagram:
$$
\xymatrix{\cdots \ar[r] &\ A_2\ar[r]\ar[d] &A_1 \ar[r]\ar[d]^{f}&A_0
\ar[r]\ar[d]
&A \ar[r]\ar[d]^{Id} &0\\
&0\ar[r] & M\ar[r] &E\ar[r] &A\ar[r] &0.}
$$
In this way  we obtain the isomorphism i).

ii) Let
$$0\to M\to C_1 \buildrel \partial \over \to C_0\to A\to 0$$
be a crossed extension. The algebra $C_0$ acts on $M$ via the
projection to $A$. Moreover
$$C_*=(\cdots \to0\to M\to C_1 \buildrel \partial \over \to C_0)$$
can be considered as a chain algebra as follows. In dimensions $0$ and $1$ it is already defined. In the dimension two one puts $C_2=M$, and  $C_i=0$
for $i>2$.  The pairing $C_i\t C_j\to
C_{i+j}$ is the given one if $i=0$ or $j=0$, while
the pairing $C_1\t C_1\to C_2$ as well as all
other pairings are zero. Then $C_*\to A$ is
an acyclic fibration. Therefore we have a lifting
$f_*:A_*\to C_*$, where $A_*\to A$ is a weak equivalence with quasi-free $A_*$.
It is clear that $f_2\in {\sf Der}(A_*,M)$ is a 2-cocycle in ${\sf Der}(A_*,M)$
and therefore gives rise
to an element in $\Sh^3(A,M)$. Conversely, starting with a 2-cocycle
$f\in {\sf Der}(A_*,M)$ one can construct the corresponding crossed extension
using the diagram
$$
\xymatrix{\cdots \ar[r]&A_3 \ar[r] \ar[d] &\ A_2\ar[r]\ar[d]^f
&A_1 \ar[r]\ar[d]&A_0
\ar[r]\ar[d]^{Id}
&A \ar[r]\ar[d]^{Id} &0\\
&0\ar[r] & M\ar[r] &C_1\ar[r]&C_0\ar[r] &A\ar[r] &0.}
$$

 The following theorem is the non-$\K$-split analogue of Theorem
 \ref{jvarediancinaag}.

\begin{The} \label{nulisgamocnoba}
The class of a crossed extension
$$0\to M\to C_1 \buildrel \partial \over \to C_0\buildrel \pi\over
\to R\to 0 \leqno(\d)$$
is zero in $\Sh
^3(R,M)$  iff $_\d\Ext(R,C_1)$ is nonempty. If this
is the case then the group $\Sh
^2(R,M)$ acts  transitively and effectively
on $_\d\Ext(R,C_1)$.
\end{The}
{\it Proof}. It is clear that the crossed extension
$$\xymatrix{ 0\ar[r]&M
\ar[r]^{\id}&M\ar[r]^{0}&R\ar[r]^{\id}&R\ar[r]&0}$$
 represents the zero element of ${\bf Cros}(R,M)$. Assume
$_\d\Ext(R,C_1)$ is nonempty and let
$$\xymatrix{&0\ar[r]&C_1\ar[d]^{\id} \ar[r]^{\mu}&S
\ar[d]^{\xi}\ar[r]^{\si}
&R\ar[d]^{\id}\ar[r]&0\\
0\ar[r]&M\ar[r]&C_1\ar[r]^\d&C_0\ar[r]^\pi&R\ar[r]&0}$$
 be an object of the category $_\d\Ext(R,C_1)$. Then $d:M\oplus
C_1\to S$ is a crossed bimodule, where $d(m,c_1)=\mu(c_1)$ and the
action of $S$ on $M\oplus C_1$ is given by
$s(m,c_1)=(\xi(s)m,\xi(s)c_1)$ and $(m,c_1)s=(m\xi(s),c_1\xi(s))$.
Then one has the following commutative diagram in ${\bf
Cross}(R,M)$:
$$\xymatrix{ 0\ar[r]&M \ar[r]^{\id}&M\ar[r]^{0}&R\ar[r]^{\id}&R\ar[r]&0\\
0\ar[r]&M\ar[d]^{\id} \ar[u]^{\id}\ar[r]^{i_1}&M\oplus
C_1\ar[d]^{p_2}\ar[u]^{p_1}\ar[r]^{\mu}&S\ar[u]^{\si}
\ar[d]^{\xi}\ar[r]^{\si}
&R\ar[d]^{\id}\ar[r]\ar[u]^{\id}&0\\
0\ar[r]&M\ar[r]&C_1\ar[r]^\d&C_0\ar[r]^\pi&R\ar[r]&0}$$ which
shows that the class of $(\d)$ in $\Sh ^3(R,M)$ is zero. Here
$p_1$ and $p_2$ are standard projections from the direct sums to
summands and $i_1$ and $i_2$ are corresponding injections.

Conversely, assume the class of $(\d)$ in $\Sh ^3(R,M)$ is zero.
 It follows from Corollary \ref{crosfiltri} that
 there exists a commutative diagram
 of crossed extensions:
$$\xymatrix{ 0\ar[r]&M \ar[r]^{\id}&M\ar[r]^{0}&R\ar[r]^{\id}&R\ar[r]&0\\
0\ar[r]&M\ar[d]^{\id} \ar[u]^{\id}\ar[r]^{i}&P_1
\ar[d]^{\ee}\ar[u]^{p}\ar[r]^{\mu}&P_0\ar[u]^{\si}
\ar[d]^{\xi}\ar[r]^{\si}
&R\ar[d]^{\id}\ar[r]\ar[u]^{\id}&0\\
0\ar[r]&M\ar[r]&C_1\ar[r]^\d&C_0\ar[r]^\pi&R\ar[r]&0.}$$ It follows
that the restriction of $\mu$ to $\Ker(p)$ is a monomorphism and
therefore we have the following commutative diagram with exact
rows:
$$\xymatrix{ & 0\ar[r]&\Ker(p)
\ar[d]^{\ee}\ar[r]^{\mu}&P_0 \ar[d]^{\xi}\ar[r]^{\si}
&R\ar[d]^{\id}\ar[r]&0\\
0\ar[r]&M\ar[r]&C_1\ar[r]^\d&C_0\ar[r]^\pi&R\ar[r]&0.}$$ One
defines the $\K$-algebra $S$ via the exact sequence
$$\xymatrix{ 0\ar[r]&\Ker(p)\ar[r]^{(-\ee,\mu)}&C_1\oplus
P_0\ar[r]&S\ar[r]&0.}$$ Here the product on $C_1\oplus P_0$ is
given by
$$(c,x)(c',x'):=(c*c'+c\xi(x')+\xi(x)c',xx').$$
One easily checks that $\Ker(p)$ is really an ideal of  $C_1\oplus
P_0$ and therefore $S$ is well-defined. Now it is clear that
$$\xymatrix{&0\ar[r]&C_1\ar[d]^{\id} \ar[r]&S
\ar[d]\ar[r]
&R\ar[d]^{\id}\ar[r]&0\\
0\ar[r]&M\ar[r]&C_1\ar[r]^\d&C_0\ar[r]^\pi&R\ar[r]&0}$$ is an
object of $_\d\Ext(R,C_1)$ and the proof is finished.

\

\

\noindent {\bf Remark}. One cannot get non-$\K$-split versions of results of
Section \ref{baumin-n}. In other words for $n>2$ neither ${\bf Extalg^n}(R,M)$ nor
${\bf Cros}^n(R,M)$ are isomorphic to $\Sh^{n+1}(R,M)$ in general. This is
because for such $n$ both groups ${\bf Extalg^n}(R,M)$ and ${\bf Cros}^n(R,M)$
vanish on injective bimodules, while Shukla cohomology does not. Indeed, if
$\K=\Z$ and $R=\F_p$, then any bimodule over $R$ is injective, while the
computation in Section \ref{shugaf}) shows that $\Sh^{2i}(\F_p/\Z,F_p)=\F_p$
for all $i$. By the same reason  the groups ${\bf Extalg^2}(R,M)$ and
$\Sh^2(R,M)$ are different.

On the other hand we have another interpretation of the higher
Shukla cohomology using chain algebras. Indeed, the above argument
can be easily modified to get the following extension of Theorem
\ref{shugaf} to higher dimensions.
%This generalization is based on the fact
%that crossed bimodules are nothing but chain algebras of length $\leq 1$.
A  chain algebra $A_*$ is called of length $\leq n$ if $A_i=0$ for
all $i>n$. Let $R$ be an algebra and $M$ be a bimodule over $R$.
For any $n\geq 1$ we let ${\bf Crosext}^{n}(R,M)$ be the category
of triples $(A_*,\alpha, \beta)$ where $A_*$ is a chain algebra of
length $\leq n$ with property $H_i(A_*)=0$ for all $0<i<n$.
Moreover $\alpha:H_0(A_*)\to R$ is an isomorphism of algebras and
$\beta:M\to H_n(A_*)$ is an isomorphism of $R$-bimodules, where
the $R$-bimodule structure on  $H_n(A_*)$ is induced via
$\alpha^{-1}$. It is clear that for $n=1$  the category ${\bf
Crosext}^{1}(R,M)$ and ${\bf Crosext}(R,M)$ are equivalent. The
argument given in the proof of part ii) of Theorem \ref{shugaf}
shows that connected components of the category ${\bf
Crosext}^{n}(R,M)$ are in one-to-one correspondence with elements
of the group $\Sh^{n+2}(R,M)$. Furthermore, for a given object
$X=(A_*,\alpha, \beta)$ of the category ${\bf Crosext}^{n}(R,M)$
one can define the category $_X\Ext(R;A_n)$ of objects
$(C_*,\gamma, \eta)$, where $C_*$ is a chain algebra of length
$\leq n$ with the property $H_i(C_*)=0$ for all $i>0$,
$\gamma:H_0(C_*)\to R$ is an isomorphism of algebras and
$\eta:C_*\to A_*$ is a chain algebra homomorphism such that the
diagram
$$\xymatrix{H_0(C_*)\ar[r]^{\gamma}\ar[d]^{\eta}&R\ar[d]^{\id}\\
H_0(A_*)\ar[r]^{\alpha}&R}
$$
commutes and $\eta_n:C_n\to A_n$ is an isomorphism. Then the category
$_X\Ext(R;A_n)$ is nonempty iff the class of $X$ in $\Sh^{n+2}(R,M)$ is zero.
If this is so, then the group $\Sh^{n+1}(R,M)$ acts transitively and
effectively on the set of components of  the category $_X\Ext(R;A_n)$.

Duskin in \cite{duskin} introduced higher torsors to obtain an interpretation
of elements of the cohomology groups in very general context. For associative
algebras his approach also gives the interpretation of $\H^3$ via crossed
bimodules, but in higher dimensions his approach is totally different from one
indicated here.

\subsection{Shukla cohomology via free crossed bimodules}

 Let $R$ be an associative algebra. We claim that there is a free crossed
module $\d:F_1\to F_0$ with ${\sf Coker}(\d)=R$. Indeed, first we take a
surjective homomorphism of rings $\pi:F_0\to R$, where $F_0$ is a free
$\K$-algebra. Then we choose a free $\K$-module $V$ together with an epimorphism
$V\to \ker(\pi)$. Finally we take $\d:F_1\to F_0$ to be the free crossed
bimodule generated by $V\to F_0$. Then $\d$ has the expected property.

\begin{Pro}
Let $R$ be an associative algebra and let $M$ be a bimodule over $R$. Let
 $$0\to E\buildrel j \over \to  F_1\buildrel \d \over \to  F_0\to R\to 0\leqno(F)$$
 be a  crossed extension with free crossed bimodule $\d:F_1\to F_0$.
Then there is an exact sequence
$$\Hom_{F_0^e}(F_1,M)\buildrel j* \over \to \Hom_{R^e}(E,M)\to \Sh^3(R,M)\to
0$$ where $j^*(h)=hj$, for $h\in \Hom_{F_0^e}(F_1,M)$.
\end{Pro}

{\it Proof}. The crossed extension $(F)$ defines an element $ e\in \H^3(R,E)$.
The homomorphism  $e_*:\Hom_{R^e}(E,M)\to \Sh^3(R,M)$ sends an element $f\in
\Hom_{R^e}(E,M)$ to $f_*(e)\in \Sh^3(R,M)$. Take any crossed extension
$$0\to M\to C_1\to C_0\to R\to 0$$
Since $F_0$ is a free algebra and $\d:F_1\to F_0$ is a free crossed bimodule,
there exists a morphism of crossed extensions
$$\xymatrix{0\ar[r] &E\ar[r]\ar[d]^f& F_1\ar[r]\ar[d]& F_0\ar[r]\ar[d]^\id&R\ar[r]\ar[d]^\id&
0\\
0\ar[r] &M\ar[r]& C\ar[r]& F_0\ar[r]&R\ar[r]& 0}
$$
which shows that $e_*:\Hom_{R^e}(E,M)\to \Sh^3(R,M)$ is an epimorphism. We
claim that  $j_*(e)=0$. Indeed, $j_*(e)$ is represented by the bottom  crossed
extension in the following diagram:
$$\xymatrix{0\ar[r] &E\ar[r]^j\ar[d]^j& F_1\ar[r]\ar[d]& F_0\ar[r]\ar[d]^\id&R\ar[r]\ar[d]^\id&
0\\
0\ar[r] &F_1\ar[r]& X\ar[r]& F_0\ar[r]&R\ar[r]& 0}
$$
Obviously $F_1\to X$ has  a retraction, hence the claim. Take any $h\in
\Hom_{F_0^e}(F_1,M)$. Then we have
$$e_*j^*(h)=(hj)_*(e)=h_*j_*(e)=0$$
Thus it remains to show that if $f\in  \Hom_{R^e}(E,M)$ satisfies $f_*(e)=0$,
then $f=hj$ for some $h\in \Hom_{F_0^e}(F_1,M)$. If $f_*(e)=0$, then we can use
Theorem \ref{nulisgamocnoba} to obtain a diagram
$$\xymatrix{0\ar[r] &E\ar[r]^j\ar[d]^f& F_1\ar[r]^\d\ar[d]^g& F_0\ar[r]^{\pi}\ar[d]^\id&
R\ar[r]\ar[d]^\id&
0\\
0\ar[r] &M\ar[r]^{j'}& C\ar[r]^{\delta}& F_0\ar[r]^{\pi}&R\ar[r]& 0\\
& 0\ar[r]& C\ar[u]^\id\ar[r]^i & S\ar[u]_t\ar[r]^p&R\ar[r]\ar[u]_\id&0}
$$
Since $F_0$ is a free $\K$-algebra, the homomorphism $t$ has a section
$s:F_0\to S$. So we have $ts=\id_{F_0}$. Since $p=\pi t$, we obtain $ps=\pi t
s=\pi$. It follows that $ps\d=\pi \d=0$, thus there exists a unique $r:F_1\to
C$ such that $s\d=ir$.  Then we have $irj=s\d j=0$ and therefore $rj=0$. On the
other hand
$$
\delta (g-r)=\delta g-\delta r=\d-tir=\d-ts\d=0.$$ Therefore there exists a
unique $h:F_1\to M$ such that $g=r+j'h$. Since $j'f=gj=rj+j'hj=j'hj$ we obtain
$f=hj$ and we are done.

\section{Some computations of Shukla cohomology}\label{shgam}
\subsection{The case $\K=\Z$}
Let $\K=\Z$ and $R=\Z/n\Z$,  $n\geq 2$. Consider the exterior
algebra $\Lambda _{\Z}^*(x)$ on a generator $x$ of degree $1$ over
$\Z$. We put $\partial (x)=n$. Then $\Lambda _{\Z}^*(x)$ is a
chain algebra, which is weakly equivalent to $\Z/n\Z$. It is clear
that the normalized Hochschild cochain complex of $\Lambda
_{\Z}^*(x)$ with coefficients in $\Z/n\Z$ has a bicomplex
structure, which is   $\Z/n\Z$ in bidegree $(i,i)$, $i\geq 0$ and
is zero elsewhere. Thus
$$\Sh^*(R/\Z,R)=R[\xi],$$
where
$$\xi\in
\Sh^2({R/\Z,R})$$ has degree $2$.
%Here and later the subscribe $k$
%indicates the fact that an object is considered over ground ring$k$.
Based on the interpretation of the second Shukla cohomology
via abelian extensions (see Section \ref{gafshu}) one easily sees
that $\xi$ represents the following extension:
$$\xi=(0\to \Z/nZ \to \Z/n^2\Z\to \Z/n\Z\to 0)\in
{\sf Extalg}_{\Z}(R,R)$$ This example can be generalized as
follows. Let $A$ be an algebra over
 $R=\Z/n\Z$. We will assume that $A$ is free as a module over $\Z/n\Z$ ( of course
this  holds automatically if $n=p$ is a prime). A ring $A_0$ is
called {\it a lifting of $A$ to $\Z$} if there exists an
isomorphism of rings $A_0/nA_0\cong A$ and additionally $A_0$ is
free as an abelian group.
\begin{Pro} Let $A$ be an algebra over
$R= {\Z/n\Z}$, which is free as an $R$-module. If $A$ has a
lifting to $\Z$ then
$$\Sh^*(A/\Z,A)\cong
  \H^*(A/R,A)[\xi]$$
\end{Pro}
{\it Proof}. Let $A_*$ be a chain algebra over $\Z$ defined  as
follows. As a graded algebra $A_*$ is the tensor product
$A_*=\Lambda _{\Z} ^*(x)\t A_0$ where $x$ has degree one. The
boundary homomorphism is defined by $\partial(x)=n$. Thus as a
chain complex $A_*$ looks as follows:
$$\cdots \to 0\to A_0 \buildrel n \over \to A_0$$
in particular $A_*\to A$ is a weak equivalence and the K\"unneth
Theorem \ref{kiuneti} for Hochschild cohomology implies
$$\Sh^*(A/\Z,A)\cong \H^*(A/R,A)\t \Sh^*(R/\Z,R)\cong
  \H^*(A/R,A)[\xi]$$

\subsection{The case $\K=\Z/p^2\Z$} Let $p$ be a prime and
$\K=\Z/p^2\Z$ and $R=\Z/p\Z$. Consider the commutative chain
algebra
$$\Lambda^*_{\Z/p^2\Z}(x)\t \Gamma^*_{\Z/p^2\Z}(y),$$ where $x$ is
of degree $1$ and  $y$ is of degree $2$. Here $\Gamma ^*$ denotes
the divided power algebra. Now we put $\partial (x)=p$ and
$\partial (y)=px$. One easily checks that in this way one obtains
a chain algebra compatible with divided powers. Since the
augmentation
$$\Lambda^*_{\Z/p^2\Z}(x)\t \Gamma^*_{\Z/p^2\Z}(y) \to \Z/p\Z$$
is a weak equivalence, one can use this chain algebra to compute
the Shukla cohomology. It is clear that
$$C^*(\Lambda^*_{\Z/p^2\Z}(x)\t \Gamma^*_{\Z/p^2\Z}(y),\Z/p\Z)\cong
C^*(\Lambda^*_{\Z/p\Z}(x)\t \Gamma^*_{\Z/p\Z}(y), \Z/p\Z)$$
where $ \Lambda^*_{\Z/p\Z}(x)\t \Gamma^*_{\Z/p\Z}(y)$ is a chain algebra with
zero boundary map. Then the K\"unneth theorem for Hochschild cohomology
 \cite{homology} implies
that
$$\Sh^*(R/\K,R)\cong R[\sigma x, \sigma y,
\sigma y^{[2]},,\cdots \sigma y^{[2^n]},\cdots],  \ \ {\rm if} \ \
p=2$$ where $|\sigma z| =1+|z|$. Similarly, if $p$
is odd, then
$$\Sh^*(R/\K,R)\cong \Lambda^*(\sigma y, \cdots
\sigma y^{[p^n]},\cdots)\t
 \Z/p\Z[\sigma x, \sigma^2 y, \cdots \sigma^2 y^{[p^n]},\cdots]\ \textrm{if $p$ is odd}$$ 
 Here we use the fact that one has an isomorphism
of algebras:
$$\Gamma _{\Z/p\Z}^*(z)\cong \Z/p\Z[z]/(z^p)\t \Z/p\Z[z]/(z^{p^2})\t \Z/p\Z[z]/(z^{p^3})\t \cdots$$
The element $\sigma x$ is still represented by the following
abelian extension of algebras
$$(0\to \Z/pZ \to \Z/p^2\Z\to \Z/p\Z\to 0)\in
{\bf Extalg}_{\Z/p^2\Z}({\Z/p\Z,\Z/p\Z})$$ while $\sigma y$ is
represented by the crossed extension of algebras:
$$(0\to \Z/pZ \to \Z/p^2\Z  \buildrel p \over \to\Z/p^2\Z\to \Z/p\Z\to 0)\in
{\bf Cros}_{\Z/p^2\Z}({\Z/p\Z,\Z/p\Z}).$$
More generally, let $A$ be an algebra over
 $\Z/p\Z$. A ring $A_0$ is called
{\it a lifting of $A$ to $\Z/p^2\Z$} if there exists an
isomorphism of algebras $A_0/pA_0\cong A$ and additionally $A_0$
is free as a $\Z/p^2\Z$-module.
\begin{Pro} Let $A$ be an algebra over
$\F_p$. If $A$ has a lifting to $\K=\Z/p^2\Z$, then
$$\Sh^*(A/\K,A)\cong
  \H^*(A/\F_p,A)\t \Sh^*(\F_p/\K,\F_p)$$
  where $$\Sh^*(\F_2/\K,\F_2)\cong \F_2[\sigma x, \sigma y,
\sigma y^{[2]},,\cdots \sigma y^{[2^n]},\cdots]$$  and
$$\Sh^*(\F_p/\K,\F_p)\cong \Lambda^*(\sigma y, \cdots
\sigma y^{[p^n]},\cdots)\t
 \F_p[\sigma x, \sigma^2 y, \cdots \sigma^2 y^{[p^n]},\cdots]$$  if
$p$  odd.
\end{Pro}
{\it Proof}. Let $A_*$ be a chain algebra over $\Z/p^2$ given as
the tensor product of chain algebras:
$$A_*=A_0\t  \Lambda^*_{\Z/p^2\Z}(x)\t \Gamma^*_{\Z/p^2\Z}(y)$$
By the K\"unneth theorem \ref{kiuneti} $A_*\to A$ is a weak
equivalence and hence
$$\Sh^*_{\Z/p^
2}(A,A)\cong \H^*_{\Z/p\Z}(A,A)\t \Sh^*_{\Z/p^2\Z}({\Z/p\Z,\Z/p\Z})$$

\

\ Let us observe that if a $\F_p$-algebra $A$ has a lifting to
$\Z$ then it has also a lifting to $\Z/p^2\Z$. It is clear that
group algebras (or more generally monoid algebras), truncated
polynomial algebras
 have
lifting to $\Z$. It is also known that any smooth commutative algebra has
lifting to $\Z$ \cite{arabia}. It is also clear that the class of algebras
having lifting to $\Z$ (or $\Z/p^2\Z$) is closed under tensor
product. It is also closed under finite cartesian products.
%There is
%a cohomological obstruction to
%the problem of lifting of a $\Z/p\Z$-algebra $R$ to $\Z/p^2\Z$
%(see Section \ref{liftobs}) which lies in $\H^3_{\Z/p\Z}(R,R)$.
%We made few computations of Shukla cohomology in case when
%algebras has lifting

\subsection{On relationship between Shukla cohomology over $\Z$
and $\Z/p^2\Z$ up to dimension three}\label{sh3} In this section
$\K=\Z/p^2\Z$ and $\H^*$ denotes the Hochschild cohomology over
$\F_p$.

Let $M$ be a bimodule over an $\F_p$-algebra $A$. Since $A$ is
also an algebra over $\Z$ and $\K=\Z/p^2\Z$, we can consider not
only the Hochschild cohomology $\H^*(A,M)$, but also the Shukla
cohomologies $\Sh^*(A/\K,M)$ and $\Sh^*(A/\Z,M)$. The ring
homomorphisms $\Z\to \K\to\F_p$ yield the natural transformations
$$\H^i(A,M)\to \Sh^i(A/\K,M)$$
and
$$b^i:\Sh^i(A/\K,M)\to \Sh^i(A/\Z,M)$$
which are obviously isomorphisms for $i=0,1$. For $i=2$, the
groups in question classify abelian extensions of $A$ by $M$,
respectively in the category of algebras over $\F_p$, $\K$ and
$\Z$. Let us observe that if $X\to Y\to Z$ is a short exact
sequence of abelian groups and $pX=0=pZ$, then $p^2Y=0$. Thus any
abelian extension of $A$ by $M$ in the category of all rings lies
in the category of algebras over $\K$. It follows that for $i=2$,
the first map $\H^2(A,M)\to \Sh^2(A/\K,M)$ is a monomorphism,
while the second homomorphism is an isomorphism:
$$b^2:\Sh^2(A/\K,M)\cong \Sh^2(A/\Z,M)$$
In higher dimensions we have
\begin{Le} For all $n$ the homomorphism $$b^n:\Sh^n(A/\K,M)\to \Sh^n(A/\Z,M)$$
is an epimorphism and it has a natural splitting.
\end{Le}

{\it Proof}. We have only to consider the case $n\geq 3$. We have
to construct the homomorphism  $d^n:\Sh^n(A/\Z,M)\to
\Sh^n(A/\K,M)$, which is a right inverse of $b^n$. We consider
more carefully the case $n=3$ and then we indicate how to modify
the argument for $n>3$. In terms of crossed extensions,
$b=b^3:\Sh^3(A/\K,M)\to \Sh^3(A/\Z,M)$ sends the class of a
crossed extension
$$0\to M\to C_1\to C_0\to A\to 0$$
of $\Z/p^2\Z$-algebras to the same crossed extension but now considered as
algebras over $\Z$.  Now we take any element from $\Sh^3(A/\Z,M)$, which is
represented by the following crossed extension of $A$ by $M$ in the category of
rings:
$$0\to M\to D_1\to D_0\to A\to 0.$$ Thanks  to Lemma \ref{crospullback} and
Corollary \ref{crosfiltri} without loss of generality one can assume that $D_0$
is free as an abelian group (this follows also from Section \ref{cmcxmod}).
Thus $V:={\sf Im}(\d)$ is also  free as an abelian group and $0\to M\to D_1\to
V\to 0$ splits as a sequence of abelian groups. It follows that $0\to M\to
D/pD\to V/pV\to 0$ is exact. On the other hand $pV$ is a two-sided ideal in
$D_0$ and therefore
 one has an exact sequence $0\to V/pV\to D_0/pV\to A$. It follows that
  $D_0/pV$ is a $\Z/p^2\Z$-algebra.
By gluing these data we get a crossed extension
$$\xymatrix{ 0\ar[r]& M\ar[r]& D_1/pD_1\ar[r]& D_0/pV\ar[r]& A\ar[r] &0}
$$
and therefore an element in $\Sh^3(A/\K,M)$. In this way we obtain the
homomorphism $d=d^3:\Sh^3(A/\Z,M)\to\Sh^3(A/\K,M)$. The commutative diagram
$$\xymatrix{0\ar[r] &M \ar[r]\ar[d]_\id& D_1\ar[r]\ar[d]&D_0\ar[r]\ar[d] &A\ar[r]\ar[d]^{\id}
& 0\\
0\ar[r]& M\ar[r]& D_1/pD_1\ar[r]& D_0/pV\ar[r]& A\ar[r] &0}
$$
shows that $bd=\id$ and the case $n=3$ is done. Assume now   $n>3$.
  According to  Remark at the and of Section \ref{gafshu} we know that elements
  of $\Sh^n(A/\Z,M)$ are equivalence classes
  of chain algebras $X_*$ of length $\leq n-2$ which are acyclic in all but the
  extreme dimensions:
  $$0\to M\to X_{n-2}\to \cdots \to X_0\to A\to 0$$
  Without loss of generality one can assume that $X_0,\cdots, X_{n-3}$ are free
as  abelian groups (use Section \ref{cmctrun}, or modify the argument in Lemma
\ref{crospullback} and Corollary \ref{crosfiltri}). By repeating the previous
argument we can construct a diagram of the form
$$\xymatrix{0\ar[r] &M \ar[r]\ar[d]_\id& X_{n-2}
\ar[r]\ar[d]& \cdots \ar[r]& X_0\ar[r]\ar[d] &A\ar[r]\ar[d]^{\id}
& 0\\
0\ar[r]& M\ar[r]& X_{n-2}/pX_{n-2}\ar[r]& \cdots \ar[r]& X_0/pV\ar[r]& A\ar[r]
&0}
$$
where $V:= \Ker(X_0\to A)$ and we are done.

 Now we analyze the kernel of the homomorphism $$b=b^3:\Sh^3(A/\K,M)\to
\Sh^3(A/\Z,M)$$

\begin{Pro} \label{langwail}
Let $A$ be an algebra over $\F_p$ and let $M$ be a bimodule over $A$. Then one
has a natural isomorphism
$$\Sh^3(A/\K,M)\cong \Sh^3(A/\Z,M)\oplus \H^0(A,M)$$
where $\K=\Z/p^2\Z$.
\end{Pro}

{\it Proof} consists of several steps. We already defined the homomorphism $
d=d^3:\Sh^3(A/\Z)\to \Sh^3(A/\K)$ with $bd=\id$. Now we define the
homomorphisms
$$e: \Sh^3(A/\K)\to \H^0(A,M), \ \  c: \H^0(A,M)\to \Sh^3(A/\K)$$ with 
$$ ed=0, \ \ ec=\id, \ \ bc=0$$
and we prove that $(b,e):\Sh^3(A/\K)\to \H^0(A,M)\oplus \Sh^3(A/\Z)$ is a
monomorphism. From these assertions the result follows.

\ {\it First step. The homomorphism $e:\Sh^3(A/\K,M)\to \H^0(A,M)$}. Let
$$0\to M\to C_1 \buildrel \partial \over \to
C_0\buildrel \pi \over \to A\to 0\leqno(\d)$$ be a crossed
extension, where $C_0$ and $C_1$ are $\K$-algebras.  Since $A$ is
an algebra over $\F_p$, one has $\pi(p1)=0$, where $1\in C_0$ is
the unit of $C_0$. Therefore one can write $p1=\d([P])$ for a
suitable $[P]$ in $C_1$. Now we put:
$$e((\d))=p[P]\in M$$
it is easy to check that $e$ is a well-defined homomorphism. Let us observe
that $e(\d)=0$ if $pC_1=0$. It follows that $ed=0$.

\ {\it Second step. The canonical class $(\sigma)_A\in \Sh^3(A/\K,A)$}.
 Let $X$ be an abelian group. We let $\Z[X]$ be the free abelian
group generated by $X$ modulo the relation $[0]=0$. Here  $[x]$ denotes an
element of $\Z[X]$ corresponding to $x\in X$.
 Then we have a canonical epimorphism $\eta:\Z[X]\to X$, $\eta([x])=0$
 which gives rise to {\it the canonical free resolution of $X$}:
$$0\to R(X)\to \Z[X]\to X\to 0$$
For any $x,y\in X$ we put
$$\brk{x,y}:=[x]+[y]-[x+y]\in R(X)$$
%and it is well-known (and easy to prove) that elements of the form
%$<x,y>$ generate $R(X)$.
We now assume that $pX=0$, that is $X$ is a vector space over
$\F_p$. By applying the functor $(-)\t \Z/p^2\Z$ to the canonical
free resolution we obtain the following exact sequence
$$0\to X \buildrel i\over \to R(X)/p^2R(X) \buildrel \sigma \over \to
 \Z/p^2\Z[X] \buildrel \eta \over \to X\to 0 \leqno(\sigma)_X$$
 Here we used the well-known isomorphism $V\cong \Tor_1(V,\Z/p^2Z)$
 for any $\F_p$-vector space $V$ considered as an abelian group
 (the $\Tor$ and $\t$ are taken of course over $\Z$ and not over
 $\K=\Z/p^2\Z$). The homomorphism $i$ has the following form
 $$i(x)=\sum_{j=0}^{p-1}p\brk{jx,x} \ \ {\sf mod} (p^2R(X))$$
Let us turn back to our situation. We can take $X=A$. The
multiplicative structure on $A$ can be extended linearly to $\Z[A]$
to get an associative algebra structure on it. Then not only
$\eta$ is a ring homomorphism, but the exact sequence $(\sigma)_A$
is a crossed extension and therefore we obtain  an element
$$(\sigma)_A=( 0\to A \to R(A)/p^2R(A) \buildrel \sigma \over \to
 \Z/p^2\Z[A]  \to A\to 0)\in \Sh^3(A/\K,A)$$
 It is clear that $A\mapsto (\sigma)_A$ is a functor from $\F_p$-algebras to the
 category of crossed extensions of $\Z/p^2\Z$-algebras. Since
 $$\sigma(\sum_{j=0}^{j=p-1}\brk{j,1})=p[1]$$
 one has $$e((\sigma)_A)=1\in \H^0(A,A)\subset A.$$ On the other hand
  $p^2R(A)$ is an ideal of $\Z[A]$. Thus  we have a commutative
diagram
$$\xymatrix{&0\ar[r] & R(A)/p^2R(A)\ar[r]\ar[d]^{\id}&
\Z[X]/p^2R(A)\ar[r]\ar[d]&A\ar[r]\ar[d]^{\id}&0\\
0\ar[r]&A\ar[r] &R(A)\ar[r] &\Z/p^2\Z[A]  \ar[r]&A\ar[r]&0 }.$$ It follows from
Theorem \ref{nulisgamocnoba} that the class $(\sigma)_A$ has the following
important property:
%In$\Sh^3(A/\Z,A)$ one holds the equality
$$0=b((\sigma)_A) \in \Sh^3(A/\Z,A).$$

\ {\bf Third step. The homomorphism $c:\H^0(A,M)\to \Sh^3(A/\K,M)$}.
 Using the class $(\sigma)_A$ we now define the homomorphism
$c:\H^0(A,M)\to \Sh^3(A/\K,M)$ by $$c(m)=f_m^*((\sigma)_A).$$ Here $m\in
\H^0(A,M)$ and $f_m:A\to M$ is the unique bimodule homomorphism with $f_m(1)=m$
and $f_m^*:\Sh^3(A/\K,A)\to \Sh^3(A/\K,M)$ is the induced homomorphism in
cohomology. Since $e$ and $b$ are natural transformations of functors it
follows that for any $m\in M$ we have
$ec(m)=ef_m^*((\sigma)_A)=f_m^*e((\sigma)_A)=f_m(1)=m$ and $ bc(m)=
bf_m^*((\sigma)_A)=f_m^*b((\sigma)_A)=0$.
 Thus
$$ec=\id  \ \ {\rm and } \ \ bc=0$$

%%%%%%%%%%%%%%%%%%%%%%%%%
\ {\it Fourth step.} It remains to show that  $$(b,e):\Sh^3(A/\K)\to
\H^0(A,M)\oplus \Sh^3(A/\Z)$$ is a monomorphism. Let
$$0\to M\to C_1\to C_0\to A\to 0$$ be a crossed
extension of $\Z/p^2\Z$-algebras which lies in $\Ker(b,e)$. Since it goes to
zero in $\Sh^3(A/\Z,M)$ one has the following diagram
$$\xymatrix{&0\ar[r]&C_1\ar[d]^{\id}
\ar[r]^{\mu}&S \ar[d]^{\xi}\ar[r]^{\si}
&R\ar[d]^{\id}\ar[r]&0\\
0\ar[r]&M\ar[r]&C_1\ar[r]^\d&C_0\ar[r]&R\ar[r]&0}$$ where $S$ is a ring. Since
$\xi$ is a homomorphisms of algebras with unit we have $[P]=p1_S$, where
$1_S$ is the unit of $S$. Therefore $e(\d)=p^21_=0$, because $(\d)$ goes also
to zero under the map $e$. It follows that $S$ is an algebra over $\Z/p^2\Z$.
Theorem \ref{nulisgamocnoba}  shows that the class of $0\to M\to C_1\to C_0\to
A\to 0$ in $\Sh^3(A/\K,M)$ is zero and the proof is finished.

%It is an algebra over $\Z/p^2\Z$
% iff $p^21_S=0$, but $p^21_S\in \H^0(R,M)\subset M$ and we are
%done

\section{A bicomplex computing Shukla cohomology}\label{relh} 
\subsection{Construction of a bicomplex}
In
this section following \cite{relshukla} we construct a
canonical bicomplex which computes the Shukla cohomology in the
special case, when the ground ring $\K$ is an algebra over a field
$k$. In this section, contrary to other parts of the paper the
tensor product $\t$ denotes $\t_k$ and not $\t _\K$. The same is
for $\Hom$.

Let $R$ be a $\K$-algebra and let $M$ be a bimodule over $R$, where $\K$ is
a commutative algebra over a field $k$. Thus $R$ is also an algebra over $k$.
We let $C^*(R,M)$ be the Hochschild cochain complex of $R$ considered as an algebra over $k$.
Similarly, we let $C^*(R/\K
,M)$ be the Hochschild cochain complex of
$R$ considered as an algebra over $\K$.
Accordingly $H^*(R,M)$ and $H^*(R/\K
,M)$ denotes
the Hochschild cohomology of $R$ with coefficients in
$M$ over $k$
 and $\K$ respectively.

We let  $K^{**}(\K,R,M)$ be the following
bicosimplicial vector space:
$$ K^{pq}(\K
,R,M)= \Hom(\K
^{\t pq}\t R^{\t q},M)$$
The $q$-th horizontal cosimplicial vector space structure comes from
the identification
$$K^{*q}(\K
,R,M)=C^*(\K
^{\t q}, C^q(R,M)),$$
where $C^q(R,M))=\Hom(R^{\t q},M)$
is considered as a bimodule over $\K
^{\t q}$ via
$$((a_1,\cdots,a_q)f(b_1,\cdots,b_q))(r_1,\cdots,r_q):=a_1\cdots a_qf(b_1r_1,\cdots,b_qr_q).$$
Here $f\in \Hom(R^q,M)$ and $a_i,b_j\in \K
$, $r_k\in R$.
The $p$-th vertical cosimplicial vector space structure comes from the
identification
$$ K^{p*}(\K
,R,M)= C^*(\K
^{\t p}\t R,M)$$
where $M$ is considered as a bimodule over $\K
^{\t p}\t R$ via
$$(a_1\t\cdots \t a_p\t r)m(b_1\t \cdots \t b_p\t s):=(a_1\cdots a_pr)m(b_1\cdots b_ps).$$
We allow ourselves to denote the corresponding bicomplex by $ K^{**}(\K
,R,M)$ as well. Thus $ K^{**}(\K
,R,M)$ looks as follows:
$$\xymatrix@C=1em{M\ar[r]^{0}\ar[d]^{\delta}& M\ar[r]^{\ Id}\ar[d]^{\delta}&M\ar[r]^{0}\ar[d]^{\delta}&\cdots\\
\Hom(R,M)\ar[r]^-d\ar[d]^{\delta}& \Hom(\K\t R,M)\ar[r]^-d\ar[d]^{\delta} &\Hom(\K\t\K\t R,M)\ar[r]\ar[d]^{\delta}&\cdots\\
\Hom(R\t R,M)\ar[r]^-d\ar[d]^{\delta} & \Hom(\K^{\t 2}\t R^{\t 2},M)\ar[r]^d\ar[d]^{\delta} &\Hom(\K^{\t 4}\t R^{\t 2},M)\ar[r]\ar[d]^{\delta}&\cdots\\
\vdots & \vdots & \vdots &}$$
Therefore for $f:\K
^{\t pq}\t R^{\t q}\to M$ the corresponding linear maps
$$d(f): \K
^{\t (p+1)q}\t R^{\t q}\to M\ \ {\rm and } \ \ \delta(f):\K
^{\t p(q+1)}\t R^{\t (q+1)}\to M$$
are given by
$$df(a_{01},...,a_{0q},a_{11},..., a_{1q},..., a_{p1},..., a_{pq},r_1,...,r_q)=$$
$$a_{01}...a_{0q}f(a_{11},..., a_{1q},..., a_{p1},..., a_{pq},r_1,...,r_q)+$$
$$+\sum_{0\leq i<p}(-1)^{i+1}f(a_{01},... ,a_{0q},...,a_{i1}a_{i+1,1},...,a_{iq}a_{i+1,q},..., a_{p1},..., a_{pq},r_1,...,r_q)+$$
$$(-1)^{p+1}f(a_{01},... ,a_{0q},...,a_{p-1,1},... ,a_{p-1,q},a_{p1}r_1,...,a_{pq}r_q).$$
and
$$\delta(f)(a_{10},...,a_{1q},..., a_{p0},...,a_{pq},r_0,...,r_q)=$$
$$(-1)^pa_{10}...a_{p0}r_0f(a_{11},...,a_{1q},a_{p1},...,a_{pq},r_1,...,r_q)+$$
$$\sum_{0\leq i<q}(-1)^{i+p+1}f(a_{10},...,a_{1i}a_{1,i+1},... ,a_{pi}a_{pi+1},..., a_{pq},r_0,...,r_ir_{i+1},...,r_q)+$$
$$(-1)^{q+p+1}f(a_{10},...,a_{1q-1},..., a_{p0},...,a_{p,q-1},r_0,...,r_{q-1})a_{1q}... a_{pq}r_q$$
We let $H^*(\K
,R,M)$ be the homology of the bicomplex  $ K^{**}(\K
,R,M)$. We also consider the following
subbicomplex $\bar{K}^{**}(\K
,R,M)$ of $ K^{**}(\K
,R,M)$:
$$\xymatrix@C=1em{M\ar[r]\ar[d]^{\delta}& 0\ar[r] \ar[d]&
0\ar[r]\ar[d]&\cdots \\
\Hom(R,M)\ar[r]^-d\ar[d]^{\delta}& \Hom(\K
\t R,M)\ar[r]^-d\ar[d]^{\delta} &\Hom(\K
\t \K
\t R,M)\ar[r]\ar[d]^{\delta}&\cdots\\
\Hom(R\t R,M)\ar[r]^-d\ar[d]^{\delta} & \Hom(\K
^{\t 2}\t R^{\t 2},M)\ar[r]^-d\ar[d]^{\delta} &\Hom(\K
^{\t 4}\t
R^{\t 2},M)\ar[r]\ar[d]^{\delta}&\cdots\\
\vdots & \vdots & \vdots &}.$$
It is clear that $H^*(\K
,R,M)\cong H^*(\bar{K}^{**}(\K
,R,M))$.

\subsection{The homomorphism $\alpha$}
It follows from the definition that
$$\Ker(d:K^{*0}\to K^{*1})\cong C^*(R/\K
,M).$$
Therefore one has the canonical homomorphism
$$\alpha^n:H^n(R/\K
,M)\to H^n(\K
,R,M), \ n\geq 0.$$

\begin{The}
 {\rm i)} The homomorphisms $\alpha^0$ and $\alpha^1$ are isomorphisms.
The homomorphism $\alpha^2$ is a monomorphism.

{\rm ii)}  If $R$ is projective over $\K
$, then
$$\alpha^n:H^n(R/\K
,M)\to H^n(\K
,R,M)$$
is an isomorphism for all $n\geq 0$

{\rm iii)} The groups $H^*(\K
,R,M)$ are canonically isomorphic to $\Sh^*(R/\K,M)$
\end{The}

{\it Proof}. i) is an immediate consequence of the definition of the bicomplex\linebreak
$\bar{K}^{**}(\K
,R,M)$. ii) The bicomplex  gives rise to the following spectral sequence:
$$E^{pq}_1=H^q(\K
^{\t p},C^p(R,M))\Longrightarrow H^{p+q}(\K
,R,M).$$
Let us recall that if $X$ and $Y$ are left modules over an associative algebra
$S$, then $\Ext_{S}^*(X,Y)\cong \H^*(S,\Hom(X,Y))$ \cite{CE},
where $\Hom(X,Y)$ is
considered as a bimodule over $S$ via $(sft)(x)=sf(tx)$.
Here $x\in X$, $s,t\in S$ and $f:X\to Y$ is a lineal map.
Having this isomorphism in mind, we can rewrite
$E^{pq}_1\cong \Ext^q_{\K
^{\t p}}(R^{\t p},M)$. By our assumptions
$R^{\t p}$ is projective over $\K
^{\t p}$. Therefore
the spectral sequence degenerates and we get $H^*(\K
,R,M)\cong
H^*(C(R/\K
,M))=H^*(R/\K
,M)$. Here we used the obvious isomorphism
$$\Hom_{\K
\t \K\t \cdots \t \K
}(R\t R\t \cdots \t R,M)=
\Hom_{\K
}(R\t_\K
 R\t_\K
 \cdots \t_\K
  R ,M).$$
iii) We let ${\bar K}^*(\K
,R,M)$ denote the total cochain complex
associated to the bicomplex
 ${\bar K}^{**}(\K
 ,R,M)$. Then this construction has an obvious extension to the category of chain $\K$-algebras.
 Unlike Lemma \ref{hhisinvariantoba}, for any weak equivalence $R_*\to S_*$ of chain $\K$-algebras
 the induced map ${\bar K}^*(\K, S_*,M) \to {\bar K}^*(\K,R_*,M)$ is  a weak equivalence. This is because
 the definition of ${\bar K}^*(\K
,R,M)$ involves the tensor products and hom's over the field $k$
and not over $\K$. Furthermore, by ii) $H^*(\K,R_*,M)$ is
isomorphic to the Hochschild cohomology, provided $R_*$ is
degreewise projective over $\K$. In particular this happens when
$R_*$ is cofibrant. Now we take any $\K$-algebra $R$ and a
cofibrant replacement $R^c_*$ of $R$. Then one has
$$H^*(\K,R,M)\cong H^*(\K,R_*^c, M)\cong \H^*(R_*^c/\K,M)\cong \Sh^*(R/\K,M).$$

\begin{Co} {\rm i)} There is a natural bijection
$${\bf
Extalg}(\K,R,M)\cong \H^2(\K
,R,M).$$

{\rm ii)} There is a natural bijection
$${\bf Cros}(\K,R,M)\cong \H^3(\K
,R,M).$$

\end{Co}
Our next aim is to describe directly the
cocycles of
 $H^*(\K,R,M)$
corresponding to abelian and crossed
extensions.

We have
$$H^2(\K
,R,M)=Z^2(\K
,R,M)/B^2(\K
,R,M)$$
where $Z^2(\K
,R,M)$  consists of pairs $(f,g)$ such that $f:R\t R\to M$
and $g:\K
\t R\to M$ are linear maps and the equalities
$$ag(b,r)-g(ab,r)+g(a,br)=0$$
$$abf(r,s)-f(ar,bs)=arg(b,s)-g(ab,rs)+g(a,r)bs$$
$$rf(s,t)-f(rs,t)+f(r,st)-f(r,s)t=0$$
hold. Here $a,b\in \K
$ and $r,s,t\in R$. Moreover, $(f,g)$ belongs to
$B^2(A,R,M)$ iff there exists a linear map $h:R\to M$ such that
$f(r,s)=rh(s)-h(rs)+h(r)s$ and $g(a,r)=ah(r)-h(ar).$
 Starting with $(f,g)\in Z^2(\K
 ,R,M)$ we construct an abelian extension of $R$ by $M$ by
putting $S=M\oplus R$ as a vector space. A $\K
$-module structure on $S$ is given
by $a(m,r)=(am+g(a,r),ar)$,
while the multiplication on $S$ is given by $(m,r)(n,s)=(ms+rn+f(r,s),rs)$. Conversely, given an abelian extension $$ 0\to M\to S \to R\to 0$$
we choose a $k$-linear section $h:R\to S$ and then we put $f(r,s):=h(r)h(s)-h(rs)$ and $
g(a,r):=ah(r)-h(ar)$. One easily checks that $(f,g)\in  Z^2(\K
,R,M)$ and one
gets i). Similarly, we have $H^3(\K
,R,M)=Z^3(\K
,R,M)/B^3(\K
,R,M)$. Here $Z^3(\K
,R,M)$  consists of triples $(f,g,h)$ such that $f:R\t R\t R\to M$,  $g: \K
\t \K
\t R\t R\to M$ and $h: \K\t \K
\t R\to M$ are linear maps and the following relations hold:
$$r_1f(r_2,r_3,r_3)-f(r_1r_2,r_3,r_4)+f(r_1,r_2r_3,r_4)-f(r_1,r_2,r_3r_4)+f(r_1,r_2,r_3)r_4=0$$
$$abc\!f(r,\!s,t)-f\!(ar\!,\!bs,ct)\!=\!arg(b,c\!,y,z)-g(ab,c,\!xy,z)+g(a,bc,x,yz)-g(a,b,x,\!y)cz$$
$$abg(c,\!d\!,x,\!y\!)-g(ac,\!bd,\!x,\!y)+g(a,b,cx,\!dy)\!=\!acxh(b,\!d,\!y)-h(ab,cd,\!xy)+h(a,c,x)bdy$$
$$ah(b,c,x)-h(ab,c,x)+h(a,bc,x)-h(a,b,cx)=0.$$
Moreover, $(f,g,h)$ belongs to
$B^3(\K
,R,M)$ iff there exist linear maps $m:R\t R\to M$ and $n:\K
\t R\to M$ such that
$$f(r,s,t)=rm(s,t)-m(rs,t)+m(r,st)-m(r,s)t$$
$$g(a,b,r,s)=abm(r,s)-m(ar,bs)-arn(b,s)+n(ab,rs)-n(a,x)bs$$
$$h(a,b,r)=an(b,r)-n(ab,r)+n(a,br).$$
%For ii) we only construct the 3-dimensional cocycle corresponding to a crossed extension.
Let $$0\to M\to C_1 \buildrel \partial \over \to C_0\buildrel \pi \over \to R\to 0$$
be a crossed extension. We put $V:={\sf Im}
(\partial)$ and consider $k$-linear sections $p:R\to C_0$
and $q:V\to C_1$
of $\pi: C_0\to R$ and $\partial:C_1\to V$ respectively. Now we define
$$m:R\t R\to V \ \ {\rm and} \ \ n:A\t R\to V$$
by $m(r,s):=q(p(r)p(s)-p(rs))$ and $n(a,r):=q(ap(r)-p(ar))$. Finally we define
$f:R\t R\t R\to M$,  $g:\K^{\t 3}
\t R\t R\to M$ and $h:\K\t \K
\t R\to M$ by
$$f(r,s,t):=p(r)m(s,t)-m(rs,t)+m(r,st)-m(r,s)p(t)$$
$$g(a,b,r,s):=p(as)n(b,s)-n(ab,rs)+bn(a,x)p(y)-abm(r.s)+m(ax,by)$$
$$h(a,b,r):= an(b,r)-n(ab,r)+n(a,bx).$$
Then $(f,g,h)\in Z^3(A,R,M)$ and the corresponding class in $\H^3(A,R,M)$
depends only on the
connected component of a given crossed extension. Thus we obtain a well-defined map
$ {\bf Cros}(A,R,M)\to \H^3(\K
,R,M)$ and a standard argument (see \cite{baumin})  shows that it is an isomorphism.

%%%%%%%%%%%%%%%%%%%%%%%%%%%%%%%%%%%%%%%%%%%%%%
\section{Applications to Mac Lane cohomology}\label{hmlsh}
%%%%%%%%%%%%%%%%%%%%%%%%%%%%%%%%%%%%%%%%%%%%%%%

In this section we are working with rings. So our ground ring  is the ring of
integers $\K=\Z$.

\subsection{Eilenberg-MacLane $Q$-construction and Mac Lane cohomology} The
definition of the Mac Lane cohomology \cite{mac} of a ring $R$ with
coefficients in an $R$-bimodule $M$ is based on the work of Eilenberg and Mac
Lane on Eilenberg-Mac Lane spaces \cite{EM}. Namely, for any abelian group $A$
Eilenberg and Mac Lane constructed a chain complex $Q_*(A)$ whose homology is
the stable homology of Eilenberg- Mac Lane spaces
$$H_q(Q_*(A))\cong H_{n+q}(K(A,n), \ \ \ n>q.$$
In low dimensions $Q_*(A)$ is defined as follows \cite{EM}, \cite{mac},
\cite{macobs}, \cite{HC}. The group $Q_0(A)=\Z[A]$ is the free abelian group
generated by elements $[a]$, $a\in A$ modulo the relation $[0]=0$. The group
$Q_1(A)$ is the free abelian group generated by pairs $[a,b]$, $a,b\in A$
modulo the relations $[a,0]=0=[0,a]$, $a\in A$, while the group $Q_2(A)$ is the
free abelian group generated by $4$-tuples $[a,b,c,d]$ modulo the relations
$$[a,b,0,0]=[0,0,c,d]=[a,0,c,0]=[0,b,0,d]=[a,0,0,d]=0$$
in general $Q_n(A)$ is generated by $2^n$-tuples modulo some relations
\cite{EM},\cite{mac},\cite{JP2}. The boundary map is given by
$$d[a,b]=[a]+[b]-[a+b]$$
$$d[a,b,c,d]=[a,b]+[c,d]-[a+c,b+d]-[a,c]-[b,d]+[a+b,c+d].$$
For any $a\in A$, the element $\gamma(a):=[0,a,a,0]\in Q_2(A)$ is a two-dimensional cycle
and $\gamma$ yields an isomorphism (see \cite{EM},\cite{mac})
$$\gamma:A/2A\cong H_2(Q_*(A)).$$

Moreover for any abelian groups $A$ and $B$ there is a natural pairing
$$Q_*(A)\otimes Q_*(B)\to Q_*(A\t B)$$
(see for example \cite{mac}, \cite{JP2} or \cite{HC}). For any ring $R$, this
pairing allows us to put a chain algebra structure on $Q_*(R)$. For example, in
very low dimensions we have
$$[x][y]=[xy],\ \ [x][y,z]=[xy,xz], \ \ [x,y][z]=[xz,yz],$$
$$[x][y,z,u,v]=[xy,xz,xu,xv], \ \ [x,y,z,t][u]=[xu,yu,zu,tu]$$
$$[x,y][u,v]=[xu,xv,yu,yv]$$

By definition the Mac Lane cohomology $\HML^*(R,M)$ is defined as the
Hochschild cohomology of $Q_*(R)$ with coefficients in $M$. One can also
introduce the dual objects -- Mac Lane homology. It was proved in \cite{PW}
that Mac Lane homology is isomorphic to the topological Hochschild homology of
B\"okstedt \cite{bo}. It is also isomorphic to  the stable $K$-theory thanks to
a result of Dundas and McCarthy \cite{dum}.

\subsection{Relation with Shukla cohomology in low dimensions}
 Since
$H_0(Q_*(R))$ $\cong$ $R$ we have a natural augmentation $\ee:Q_*(R)\to R$. Since
$Q_*(R)$ is free as an abelian group the chain algebra
$$V_*(R)=(\cdots \to 0 \to \ker (\ee)\to Q_0(R))$$
is $\Z$-free and $V_*(R)\to R$ is a weak equivalence. Hence $V_*(R)$
can be used to
compute the Shukla cohomology. Thus the morphism of chain algebras
$$
\xymatrix{\cdots \ar[r]&Q_2(R)\ar[r]\ar[d]&Q_1(R)\ar[r]\ar[d]& Q_0(R)\ar[r]\ar[d]^{\id}
&R\ar[r]\ar[d]^{\id}&0\\
&0\ar[r]& \Ker(\ee)\ar[r]&Q_0(R)\ar[r]&R\ar[r]&0}
$$
yields the natural transformation
\begin{equation}\label{shumac}
\Sh^i(R/\Z,M)\to \HML^i(R,M)
\end{equation}
which is an isomorphism in dimensions 0,1 and 2. Thus $\HML^2(R,M)$ classifies
singular extensions of $R$ by $M$ in the category of rings, see also
\cite{mac}. According to Theorem 9 of \cite{macobs} in the dimension $3$ one
has the following exact sequence (see also Theorem \ref{zustimimdevrobebi})
\begin{equation}\label{HMLorigr}
0\to \Sh^3(R/\Z,M)\to \HML^3(R,M)\to \H^0(R, \ _2M)
\end{equation}
The connecting map $\HML^3(R,M)\to \H^0(R, \ _2M)$ is defined via $\gamma$ (see
\cite{macobs}).
\begin{Pro}\label{HMLpalg}
Let
 $R$ be an algebra over $\F_p$
  and $M$ be an $R$-bimodule. Then the natural map
$$\Sh^3(R/\Z
,M)\to \HML^3(R,M)$$ is an isomorphism.
\end{Pro}

{\it Proof}. If $p\not =2$ then this
 is an immediate consequence of the exact sequence (\ref{HMLorigr}),
 because $_2M=0$. So we have to consider only the case $p=2$. For any
$\F_2$-algebra $R$ we have the canonical homomorphism $\F_2\to R$, which yields
the following commutative diagram
$$\xymatrix{0\ar[r]&\Sh^3(R/\Z,M)\ar[r]\ar[d]&\HML^3(R,M)\ar[r]\ar[d]&\H^0(R,M)\ar[d]
\\
0\ar[r]&\Sh^3(\F_2/\Z ,M)\ar[r]&\HML^3(\F_2,M)\ar[r]&\H^0(\F_2,M)=M }
$$
It is well known that $\HML^3(\F_2,M)=0$ see for example \cite{fls} or \cite{bo}.
 Since the last vertical arrow is a monomorphism
 we are done.

Based on Proposition \ref{HMLpalg} and Proposition \ref{langwail} we obtain the
following
 \begin{Co}\label{hjb}
 Let $A$ be an algebra over $\F_p$ and let
$M$ be an $A$-bimodule.
 Then one has a split exact sequence
$$0\to \H^0(A,M)\to \Sh^3(A/\K,M)\to \HML^3(A,M)\to 0$$ where $\K=\Z/p^2\Z$.
\end{Co}
{\bf Remark}. The homomorphism $\Sh^3(R/\Z ,M)\to \HML^3(R,M)$ in general is
not an isomorphism. For example, if $R=\Z$, then $\Sh^i(\Z,-)=0$ for all $i\geq
1$, thanks to Lemma \ref{hhisinvariantoba}. On the other hand $\HML^*(\Z,-)$ is
quite nontrivial (see \cite{bo}, \cite{FP}) and in particular
$\HML^3(\Z,\F_2)=\F_2$. More about $\HML^*(\Z,-)$ see  at the end of Section
\ref{ssshmac}.

\subsection{Relation with Shukla cohomology in higher dimensions}\label{ssshmac}  The
relationship between Shukla  cohomology $\Sh^*(A/\Z,M)$ and Mac Lane cohomology
$\HML^*(A,M)$ in higher dimensions is more complicated. Let us first consider
the crucial case $A=\Z/p^k\Z$. We already saw   $\Sh^i(A/\Z,M)$, $A=\Z/p^k\Z$,
is $M$ if $i$ is even and is zero otherwise. Unlike the Shukla cohomology,
the behavior of $\HML^*(A,M)$ depends on whether $k=1$ or $k>1$. If $k=1$, then
similarly to Shukla cohomology the group $\HML^i(A,M)$ is $M$ if $i$ is even
and is zero otherwise. However the natural map
$$\Sh^i(\F_p/\Z, M)\to \HML^i(\F_p, M)$$
is an isomorphism only for $i=0, \cdots, 2p-1,$ and it is zero for
 $i>2p-2$. This
follows from the fact that $\Sh^*(\F_p, \F_p)$ is a polynomial algebra on the
generator $x$  of dimension two and $\HML^*(\F_p, \F_p)$ is a divided power
algebra on the same generator $x$ \cite{fls}. If $k>1$, then situation with Mac
Lane cohomology is more complicated. A computation made in \cite{P4} shows that
$$\HML^{2n}({\Z/p^k\Z}, \F_p)= (\F_p)^{t}, \ \ \ \HML^{2n-1}({\Z/p^k\Z},\F_p)=(\F_p)^s,$$
where $t=1+[\frac{n}{p}]$ and $s=[\frac{n+1}{p}]$. The full computation of
$\HML_*({\Z/p^k\Z},-)$ was obtained  by Brun \cite{brun}.

The relationship between Mac Lane cohomology and Shukla cohomology  for general
rings in all dimensions is given by the following theorem proved in \cite{PW}
(see also \cite{P3}).

\begin{The}\label{zustimimdevrobebi}  Let $R$ be a ring.
Then for any $R$-bimodule $M$ there is a  spectral sequence
$$E^{pq}_2(\K)=\Sh^p(R/\K,\HML^q(\K,M))\Longrightarrow \HML^{p+q}(R,M)$$
which is natural in $R$ and $M$. The spectral sequence in low dimensions gives
rise to the exact sequence:\def\sto{{\scriptstyle\to}}
$$0\sto\Sh^3(R/\Z,\!M)\sto\!\HML^3(R,\!M)\sto\H^0(R,_2\!M)\sto\!\Sh^4(R/\Z,M)\sto\!\HML^4(R,M).$$
%iii) Let
 %$A$ be an algebra over $\F_p$ and $M$ be a $R$-bimodule. Then there is an exact sequence
%$$ 0\to \H^2(A,M)\to \HML^2(A,M) \buildrel \partial \over \to \H^0(A,M)\buildrel d \over \to
%\H^3(A,M)\to \HML^3(A,M)\to $$
%$$\to \H^1(A,M)\to \H^4(A,M)\to \HML^4(A,M).$$
\end{The}

For the proof of the first part  we refer to \cite{PW} and \cite{P3}. The second
part first was proved in  \cite{JP1}. It is an immediate consequence of the
existence of the spectral sequence together with the following computation due
to B\"okstedt \cite{bo} (see also \cite{fls} and \cite{FP} ).
$$\HML^{2n}({\Z},M)= M/nM, \ \ \ \HML^{2n-1}({\Z},M)= \ _nM, \ n>0.$$
$$\HML^{2n}({\F_p}, M)= M, \ \ \ \HML^{2n-1}({\F_p},M)=0.$$

\subsection{Mac Lane cohomology and cohomology of small categories}\label{bwjp}
In this section we recall the relationship between  Mac Lane cohomology and
cohomology of small categories  \cite{JP2}. We assume that the reader is familiar with
definition of cohomology of small categories with coefficients in a natural system \cite{BW},
\cite{BJP}. Let us recall that any bifunctor gives rise to a natural system, and therefore we
can talk about the cohomology of small categories with coefficients in a bifunctor.

For a ring $R$ we let $R\mod$ be
the category of  finitely generated free $R$-modules. Actually we will assume
that objects of $R\mod$ are natural numbers and morphisms from $\bf n$ to $\bf
m$ are the same as $R$-linear maps $R^n\to R^m$, or $m\times n$-matrices over
$R$. Let $M$ be a bimodule over $R$. There is a bifunctor
$$\hom(-,M\t_R-):R\mod^{op}\times R\mod\to \sf Ab$$ given by
$$\hom(-,M\t_R-)(X,Y)=\Hom_R(X,M\t _RY).$$
Therefore one can consider the cohomology $H^*(R\mod, \hom(-,M\t_R-))$ of the
category $R\mod$ with coefficients $\hom(-,M\t_R-)$ in the sense of Baues and
Wirshing \cite{BW} (see also \cite{BJP}). A result of \cite{JP2} asserts that
one has an isomorphism:
\begin{equation}\label{jibpira}
\HML^*(R,M)\cong H^*(R\mod, \hom(-,M\t_R-)).
\end{equation}
Comparing this isomorphism with the natural homomorphism $\Sh^*(R,M)\to
\HML^*(R,M)$ one obtains the homomorphism
\begin{equation}\label{shbw}
\Sh^i(R,M)\to H^i(R\mod, \hom(-,M\t_R-)), \ \ i\geq 0.
\end{equation}
Now we recall the description of this homomorphism in terms of extensions for
$i=2.$
 Let
$$0 \to M\buildrel i\over\to S\buildrel p\over \to R\to 0$$
be an abelian extension of rings. Then
$$0\to \hom(-,M\t_R-)\to S\mod \buildrel p_*\over \to R\mod \to 0$$
is a linear extension of categories \cite{BW, BJP}, where  the functor
$p_*$ is given by $p_*(A)= A\t_SR$, $A\in S\mod$ (having in mind the
identification of $R\mod$ as the category of natural numbers and matrices, the
functor $p_*$ is the identity on objects and is given by applying $p$ on
matrices). Let us recall that for fixed $R$ and $M$ the equivalence classes of
abelian extensions of $R$ by $M$ form a group ${\sf Extal}(R,M)$, which is
isomorphic to the second Shukla cohomology of $R$ with coefficients in $M$ (see
Theorem \ref{shugaf}), while linear extensions are classified using the second
cohomology of small categories \cite{BW},\cite{BJP}, thus we obtain the
homomorphism
$$\Sh^2(R,M)\to \H^2(\mod(R),\hom(-,M\t_R-))$$
which is an isomorphism according to isomorphisms (\ref{jibpira}) and
(\ref{shumac}).  One easily shows that any biadditive bifunctor $D$ on $R\mod$
is of the form $D=\hom(-,M\t_R-)$, where $M=D(R,R)$. Thus one can conclude that
any extension of the category $R\mod$ by a biadditive bifunctor is also of the
form $S\mod$, for some ring $S$. In particular it is an additive category, more
generally any linear extension of an additive category by biadditive functor is
an additive category. This fact is an immediate consequence of
% Lemma \ref{compatible}.
Lemma 5.1.2 of \cite{BJP}.
%%%%%%%%%%%%%%%%%%%%%%%%%%%%%%%%%%%%%%%%%%%%%%%%%%%%%%%%%%%%%%%%%%%%%%%%%%
%%%%%%%%%%%%%%%%%%%%%%%%%%%%%%%%%%%%%%%%%%%%%%%%%%%%%%%%%%%%%%%%%%%%%%%%%
%%%%%%%%%%%%%%%%%%%%%%%%%%%%%%%%%%%%%%%%%%%%%%%%%%
\section{Applications to  strengthening of
additive track theories}\label{shstrict}
%%%%%%%%%%%%%%%%%%%%%%%%%%%%%%%%%%%%%%%%%%%%%%%%%%%%%%
%%%%%%%%%%%%%%%%%%%%%%%%%%%%%%%%%%%%%%%%%%%%%%%%%%%%%%%%%%%%%%%%%%%%%%%%%%%%%
%%%%%%%%%%%%%%%%%%%%%%%%%%%%%%%%%%%%%%%%%%%%%%%%%%%%%%%%%%%%%%%%%%%%%%%%%%%%%%
\subsection{Additive and very strongly additive
track theories}\label{ota} Let us recall that a \emph{track category} $\ta$ is
a category enriched in groupoids \cite{BJP}. Thus $\ta$ consists of objects and
for each pair of objects $X$, $Y$ of $\ta$ there is given the \emph{Hom-groupoid}
$\hog{X,Y}$, whose objects are termed maps, while 2-arrows --- tracks. To any
track category $\ta$ there is an associated category $\ta_\ho$ with the same
objects as $\ta$, while for objects $A$ and $B$  of $\Ob(\ta)$
 the set of morphisms $[A,B]$  in $\ta_\ho$ is the set of
connected components of the groupoid $\hog{X,Y}$.

 A track category is \emph{abelian} if for any
1-arrow $f:X\to Y$, the group $\Aut(f)$ of tracks from $f$ to itself is
abelian.  Any abelian  track category defines a natural system $D=D_{\ta}$ on
$\ta_\ho$ and a canonical class $\Ch(\ta)\in H^3(\ta_\ho,D)$ --- see section 2.3 of
\cite{BJP}. Conversely for any category $\bf C$, any natural system $D$ on $\bf
C$ and any element $a\in H^3({\bf C},D)$ there exists an abelian track category
$\ta=\ta_{{\bf C},D,a}$ unique up to equivalence such that $\ta_\ho=\bf C$ and
$\Ch(\ta) =a$ (see \cite{BJP}). In fact for a given natural system $D$ on a
category $\bf C$ there is a category ${\sf Trext}({\bf C},D)$ whose objects are
abelian track categories $\ta$ with $\ta_\ho=\bf C$ and $D_{\ta}=D$ and the set
of connected components of ${\sf Trext}({\bf C},D)$ is isomorphic to the third
dimensional cohomology \cite{P1},\cite{P2}:
\begin{equation}\label{mompares}
\pi_0({\sf Trext}({\bf C},D))\cong H^3({\bf C},D)
\end{equation}

 A {\it lax coproduct} $A\vee B$ in a track category $\ta$ is an object
$A\vee B$ equipped with maps $i_1:A\to A\vee B$, $i_2:B\to A\vee B$ such that
the induced functor
$$(i_1^*,i_2^*):\hog{A\vee B,X}\to \hog{A,X}\x \hog{B,X}$$
is an equivalence of groupoids for all objects $X\in \ta$. The coproduct is
{\it strong} if the functor $(i_1^*,i_2^*)$ is an isomorphism of groupoids.  By
duality we have also the notion of
 {\it lax product} and  {\it strong product}.

A {\it lax zero object}
 in a track category $\ta$ is an object $0$ such that the categories
 $\hog{0,X}$ and $\hog{X,0}$ are equivalent to
 the trivial groupoid for all
 $X\in\ta$. Let us recall that a \emph{trivial groupoid}
 has only one object and one arrow. A {\it strong zero object}
 in a track category $\ta$ is an object $0$ such that the categories
 $\hog{0,X}$ and $\hog{X,0}$ are trivial groupoids.

A \emph{theory} is a category possessing finite products. A \emph{track theory}
(resp. {\it strong track theory}) is a track category $\ta$ possessing finite
lax products (resp. strong products) \cite{BJP}. If $\ta$ is a track theory,
then $\ta_\ho$ is a theory. In this case the corresponding natural system on
$\ta_\ho$ is a so called \emph{cartesian natural system}, meaning that it is
compatible with finite product in an appropriate sense \cite{BJP}. Conversely,
if $\ta$ is a track category, with property that $\ta_\ho$ is a theory and
corresponding natural system is  a cartesian natural system then $\ta$ is a
track theory \cite{BJP}.

Morphisms of track theories are enriched functors which are compatible with lax
products. An equivalence of track theories is a track theory morphism which is
a weak equivalence \cite{BJP} and two track theories are called
\emph{equivalent} if they are made so by the smallest equivalence relation
generated by these. Two track theories $\ta$ and $\ta'$ are equivalent iff
there is an equivalence of categories $\ta_{\ho}\cong
\ta'_\ho$ and after identification of these categories one should have $D_{\ta}=D_{\ta'}$ and
${\sf Ch}(\ta)={\sf Ch}(\ta')$.

 The main result of \cite{BJP} is the so called \emph{strengthening theorem}, which asserts that
any abelian track theory is equivalent to a strong one.

An \emph{additive track theory} is a track category $\ta$ such that $\ta_\ho$ is an additive
category and the corresponding natural system is a biadditive bifunctor.
We are going now to give an equivalent definition, but let us before that discuss the definition of an
additive category. Let $\A$ be a category with zero object $0$ which possesses
also finite coproducts and finite products. For objects $A$ and $B$ we have
canonical inclusions $i_1=(\id,0):A\to A\x B$ and $i_2=(0,\id):B\to A\x B$  and
therefore also the canonical morphism $\kappa:A\vee B\to A\x B$.
 The category $\A$ is called \emph{semi-additive}
if the canonical morphism $\kappa:A\vee B\to A\x B$  is an isomorphism for all
$A$ and $B$. If $\A$ is a semi-additive category and $f,g:A\to B$ are morphisms
in $\A$, we let $f+g:A\to B$ be the following composite:
$$\xymatrix{A\ar[r]^{\Delta}&A\x A\ar[r]^{(f,g)}&B\x
B\ar[r]^{\kappa^{-1}}&B \vee B\ar[r]^\nabla&B}$$ where $\Delta=(\id,\id)$ is
the diagonal and $\nabla$ is the codiagonal. Thus in a semi-additive category
hom's are commutative monoids and the composition law is  biadditive. If
these monoids are abelian groups then a semi-additive category is called
\emph{additive}. This happens iff the identity morphism $\id_A$ admits the
additive inverse $-\id_A$, for each object $A$.

Now we pass to the 2-world. Let $\ta$ be a track theory with lax zero object. Then
for any objects $A$ and $B$ of $\ta$, there is a map $i_1:A\to A\x B$ and
tracks $p_1i_1\then \id_A$, $p_2i_1\then 0$. Similar meaning has $i_2:B\to A\x
B$. A {\it semi-additive track category} is an additive track theory with strong
zero object, such that for any two objects $A$ and $B$ the lax product $A\x B$
is also lax coproduct via $i_1:A\to A\x B$ and $i_2:A\to A\x B$. It is clear
that the homotopy category $\ta_\ho$ of a semi-additive track theory is a
semi-additive category.

One can prove that a track category $\ta$ is an additive track theory iff it is
a semi-additive track category and additionally the semi-additive category
 $\ta_{\ho}$ is an additive category.

An additive track category is called {\it very strong} if it admits strong zero
object $0$, strong finite products and for any two objects $A$ and $B$ the
strong product $A\x B$ is also the strong coproduct by $i_1:A\to A\x B$ and
$i_2:A\to A\x B$.

As we said a strengthening theorem of \cite{BJP} asserts that
any track theory is equivalent to a strong one. In particular, any
additive track category is equivalent to one which possesses strong products.
Since the dual of an additive track category is still a track theory, we see
that it is also equivalent to one which possesses strong coproducts. Can we always
get strong products and coproducts simultaneously? In other words, is every
additive track category $\ta$ equivalent to a very strong one? We will see
that the answer is negative in general, but positive provided the corresponding
homotopy category $\ta_\ho$ (which is an additive category in general) is
$\F_2$-linear, or $2$ is invertible in $\ta_\ho$ (meaning that all $\Hom$'s are
modules over  $\Z[\frac{1}{2}])$. More precisely the following is true:

\begin{The} \label{verystongstrengthening} Let $\ta$ be a small
additive track theory with the homotopy category ${\bf C}=\ta_\ho$ and a
canonical bifunctor $D=D_{\ta}$. Let $_2D$ be the two-torsion part of $D$. Then
there is a well-defined element $o(\ta)\in H^0({\bf C},\ _2D)$, which is
nontrivial in general and such that $o(\ta)=0$ iff $\ta$ is equivalent to a
very strongly additive track theory. The class $o(\ta)$ is zero provided hom's of
the additive category ${\bf C}$ are modules either over $\Z[\frac{1}{2}])$ or over $\F_2$.
\end{The}

The reader should compare Theorem \ref{verystongstrengthening} with the exact
sequence (\ref{HMLorigr}) and Proposition \ref{HMLpalg}. The similarity of
these results is not accidental. Indeed, let us give a quick proof of the Theorem
\ref{verystongstrengthening} in the the key case when ${\bf C}=R\mod$ is the
category  of finitely generated free modules over a ring $R$.

The proof of Theorem \ref{verystongstrengthening} in the general case is a
repetition of the proof given below in the special case, except that one has to use ringoids
instead of rings and we leave it as an exercise to the interested reader.

{\it Proof of  Theorem \ref{verystongstrengthening}. The case ${\bf C}=R\mod$}. For any
biadditive bifunctor $D$ on $R\mod$ one has an isomorphism $D\cong
\hom(-,M\t_R-)$, where $M=D(R,R)$. Here  we used the notations of Section
\ref{bwjp}. By Proposition \ref{verysh} $\Sh^3(R,M)$ classifies all very strong
additive track categories (up to  equivalence) $\ta$ with $\ta_\ho=R\mod$ and
$D(R,R)=M$, where $D$ is the canonical bifunctor associated with $M$. On the
other hand the isomorphism (\ref{mompares}) and the isomorphism (\ref{jibpira})
show that $\HML^3(R,M)$ classifies all additive track categories $\ta$ (up to
equivalence) with $\ta_\ho=R\mod$ and $D(R,R)=M$. Let $\ta$ be an additive
track category, then up to isomorphism (\ref{jibpira}) one  can assume that
${\sf Ch}(\ta)\in \HML^3(R,M)$. Thanks to the exact sequence (\ref{HMLorigr})
we can take $o(\ta)$ to be the image of ${\sf Ch}(\ta)$ in $\H^0(R, \ _2M)$.
Now Theorem \ref{verystongstrengthening} is  a consequence of the exact
sequence (\ref{HMLorigr}) and Proposition \ref{HMLpalg}. The example $R=\Z$ and
$M=\F_2$ shows that the map $\HML^3(R,M)\to \H^0(R,M)$ is not trivial in
general. It follows that the function $o$ is not trivial in general.

{\bf Remarks}. 1) The following example introduces a well-known example from
topology \cite{B1} of a track category $\ta$, which represents the generator of
$\HML^3(\Z,\F_2)=\F_2$. Following \cite{B1} we consider the track category
${\sf Top}^*$ of compactly generated Hausdorff spaces with basepoint $*$. Maps
in ${\sf Top}^*$ are pointed maps. A track $\aa:f\then g$ between pointed
maps $f,g:A\to B$ is a homotopy class of a homotopy relative to $A\times
\partial I$. Now we take ${\cal S}^k$ to be the full subcategory of ${\sf Top}^*$ consisting
of finite one-point unions of spheres $S^k$, $k\geq 2$. Then ${\cal S}^k$ is
an abelian track category and $({\cal S}^k)_\ho$ is equivalent to $\Z\mod$. For
$k\geq 3$ the corresponding bifunctor is $\hom(-,\F_2\t -)$ and therefore
 ${\cal S}^k$ is an additive track theory, whose class in
 $$H^3(\Z\mod, \hom(-,\F_2\t -))\cong \HML^3(\Z,\F_2)=\F_2$$
 is nontrivial.

\noindent  2) One can describe the function $o$ in Theorem
\ref{verystongstrengthening} as follows. Let $\ta$ be an additive track theory.
Let $\vee$ denote the weak coproduct in $\ta$ and let $0$ be the weak zero object. For
objects $X,Y$ one has therefore ``inclusions'' $i_1:X \to X\vee Y$ and $i_2:Y \to
X\vee Y$. Since $X\vee Y$ is also a weak product of $X$ and $Y$ in $\ta$ it
follows that one has also projection maps $p_1:X\vee Y\to X$ and $p_2:X\vee
Y\to Y$.
%For any $X$ we choice maps $0_X:X\to 0$ and $0^X:0\to X$. We also
For each $X$ we choose maps $i_X:X\to X\vee X$ and $t_X:X\vee Y\to Y\vee X$ in
such a way that classes of $i_X$ and $t_X$ in  $\ta_\ho$ are the codiagonal and
twisting maps in the additive category $\ta_\ho$. It follows that there is a
unique track
$$\aa_X:i_X\then t\circ i_X$$
such that $p_{i*}(\aa_X)=0$ for $i=1,2$. Now, let $(1,1):X\vee X\to X$ be a map
which lifts the codiagonal map in $\ta_\ho$. Then $(1,1)_*\aa_X$ is a track
$\id_X\to \id_X$ and therefore it differs from the trivial track by an element
$o(X)\in D(X,X)$. One can prove that the assignment $X\mapsto o(X)$ is an
expected one.

3) Corollary \ref{HMLpalg} shows that if $\ta$ is an additive track theory such that
$\ta_\ho$ is an $\F_p$-linear category then $\ta$ is equivalent to a strong
additive track theory with $\Z/p^2\Z$-linear hom's. This fact for a special
track theory arising in the theory of the ``secondary Steenrod algebra'' was proved
by the first author by completely different methods and was a starting point of
this work.

4) Based on quadratic categories and square rings \cite{square} we in the
fortcomming paper we introduce the notion of strongly additive track theories
and we prove that any additive track category is equivalent to srong one.

\subsection{Crossed bimodules and very strongly additive track theories}
Let us recall that for a category $\bf C$ and a  bifunctor $D$ there is a
category ${\sf Trext}({\bf C},D)$ such that $\pi_0({\sf Trext}({\bf C},D))\cong
H^3({\bf C},D)$. The objects of ${\sf Trext}({\bf C},D)$ are abelian track
categories $\ta$ with $\ta_\ho\cong \sf C$ and $D_{\ta}=D$. If additionally
$\bf C$ is an additive category and $D$ is a biadditive bifunctor, then any
such $\ta$ is an additive track theory. We let ${\sf Strext}({\bf C},D)$ be the
full subcategory of ${\sf Trext}({\bf C},D)$ whose objects are very 
strongly additive
track theories.

\begin{Pro}\label{verysh} Let $R$ be a ring and let $M$ be a bimodule over $R$.
There is a functor
$${\bf Crosext}(R,M)\to {\sf Strext}(R\mod,\hom(-,M\t _R-))$$
which is an equivalence of categories.
\end{Pro}

{\it Proof}.  Let
$\d:C_1 \to C_0$ be a crossed bimodule. We let $\ta=\ta(\d)$ be the
following track category. The objects of $\ta$ are
the same as the objects of $R\mod$,
i.~e. natural numbers. For any natural numbers $n$ and $m$
the maps from $n$ to $m$
(which is the same as objects of
the groupoid $\ta(n,m)$)  are $m\times n$-matrices with coefficients in $C_0$. For
$f,g\in Mat_{m\times n}(S)$ the set of tracks $f\to g$ (which is the same as the set of
 morphisms from $f$ to $g$ in the groupoid $\ta(n,m)$) is given by
$$\Hom_{\ta(n,m)}(f,g)=\left\{h \in Mat_{m\times n}(C_1)\ \mid\ \d (h)=f-g\right\}.$$
The composition of 1-arrows is given by
the usual multiplication of matrices, while the
composition of tracks is given by the addition of matrices. One easily checks
that in this way one really  obtains a very strongly additive track theory
$\ta(\d)$. It is clear that $\ta_\ho=R\mod$, where $R={\sf Coker}(\d)$ and the
bifunctor associated to $\ta$ is $D=\hom(-,M\t_R-)$. Thus we obtain
a functor $${\bf Crosext}(A,M)\to {\sf Strext}(R\mod,\hom(-,M\t _R-)).$$
Now we construct the functor in the opposite direction.
%\begin{Le}
%Let $p:S\to R$ be a surjective homomorphism, and let $M$ be an
%$R$-$R$-bimodule. Then the restriction of the above map yields the bijection
%$${\sc Cros}(R,S,;M)\to \H^3(\mod(R),\mod(S);  \hom(-,M\t_R-)).$$
%\end{Le}
%%%%%{\it Proof}.
Let $\ta$ be an object of ${\sf Strext}(R\mod,\hom(-,M\t _R-))$. Let $\ta_0$ be
the category with the same objects as $\ta$ and with maps (i.e. 1-arrows) of
$\ta$ as morphisms. Since $\ta$ is a strongly additive track theory, we see
that $\ta_0$ is an additive category and therefore it is equivalent to $S\mod$,
where $S={\sf End}_{\ta_0}(1)$. The restriction of the quotient functor $\ta\to
\ta_0$ yields the homomorphism of rings $S\to R$. One defines $X$ to be the set
of pairs $(h,x)$, where $x\in \Hom_{\ta_0}(1.1)$ and $h:x\then 0$ is a track in
the groupoid $\ta(1,1)$. Moreover we put $\partial=\partial_{\ta}(h,x)=x$.
Then $X$ carries a structure of a bimodule over $S$, and $$0\to M\to X\buildrel
\partial \over \to S\to R\to 0$$ is a crossed extension. Then $\ta\mapsto
\partial_{\ta}$ yields the functor
$$
\mathsf{Strext}(R\mod,\hom(-,M\t _R-))\to\mathbf{Crosext}(A,M).
$$
 One easily checks that these two functors
yield the expected equivalence of categories.

%%%%%%%%%%%%%%%%%%%%%%%%%%%%%%%%%%%%%%%%%%%

\appendix

%%%%%%%%%%%%%%%%%%%%%%%%%%%%%%%%%%%%%%%%%%%%%%%%%%%%%%%%%%%%%%%%%%%%%%%%

\section{Closed model category structure on chain algebras and crossed
bimodules }\label{cmc}
%%%%%%%%%%%%%%%%%%%%%%%%%%%%%%%%%%%%%%%%%%%%%%%%%%%%%%%%%%%%%%%%%%%%%%%
%%%%%%%%%%%%%%%%%%%%%%%%%%%%%%%%%%%%%%%%%%%%%%%%%%%%%%%%%%%%%%%%%%%%%%%%
\subsection{Closed model categories}
%%%%%%%%%%%%%%%%%%%%%%%%%%%%%%%%%%%%%%%%%%%%%%%%%%%%%%%%%%%%%%%%%%%%%%%%W
We recall the definition of a closed model
 category introduced by Quillen \cite{qu}. We refer the reader
 to \cite{dwyerspalinski}  for the basic facts on the closed model
category theory.
  Let $\c$ be a category. A
morphism $f$ is a \emph{retract} of a morphism $g$ if there exists
a commutative diagram of the form
$$
\xymatrix{A\ar[r]\ar[d]_f &C\ar[d]^g\ar[r]& A\ar[d]^f \\
 B\ar[r] &  D\ar[r] &A}$$
in which the horizontal composites are identities. Let $i:A\to B$ and
$p:X\to Y$ be morphisms in $\c$. Then $i$ has \emph{left lifting
property} with respect to $p$ and $p$ has \emph{right lifting
property} with respect to $i$, if for every commutative diagram
$$
\xymatrix{A\ar[r]^{ g'}\ar[d]_i &X\ar[d]^f \\
 B\ar[r]_{ g} &  Y
}$$ there exists a commutative diagram
$$
\xymatrix{A\ar[r]^{ g'}\ar[d]_i &X\ar[d]^f \\
 B\ar[r]_{ g}\ar[ru]_h &  Y
}$$ Then $h$ is called a {\it lifting}.

\begin{De} A  closed model category consists
of  a category  ${\cal C}$ together with three distinguished
classes of morphisms  called respectively weak equivalences,
cofibrations and fibrations, so that the following  5 axioms hold.

\smallskip
\noindent {\bf CM 1}. $\cal C$ has all finite limits and colimits. All
 3 classes form a subcategory.

\smallskip
\noindent {\bf CM 2}. If $f$ and $g$ are composable arrows in
$\cal C$ and two of the three morphisms $f,g,gf$ are weak
equivalences, then so is the third.

\smallskip
\noindent {\bf CM 3}. A retract of a fibration (resp. cofibration,
weak equivalence) is still a fibration (resp. cofibration, weak
equivalence).

\smallskip
\noindent {\bf CM 4}. Fibrations have the right lifting property
with respect to acyclic cofibrations and cofibrations have left
lifting property with respect to acyclic fibrations. Here a map is
called an \emph{acyclic fibration} (resp. \emph{acyclic
cofibration}) if it is both a fibration (resp. cofibration) and a
weak equivalence.

\smallskip
\noindent {\bf CM 5}. Any arrow $f:A\to B$ has factorizations
$f=pi$ and $f=qj$, where $i$ and $j$ are cofibrations, $p$ and $q$
are fibrations and $p$ and $j$ are weak equivalences too.
\end{De}

 Here is more language corresponding to closed model categories.
 An object $X$ is called {\it cofibrant} if $
\emptyset\to X$ is a cofibration. An object $Y$ is called {\it
fibrant} if $ Y\to *$ is a fibration. Here $ \emptyset$ and $*$
are respectively initial and terminal objects in $\cal C$. For any
object $X$ there are weak equivalences $X\to X^f$ and $X^c\to X$
with fibrant $X^f$ and cofibrant $X^c$. This is an easy
consequence of CM 5. Any such $X^c$ (resp. $X^f$) is called a
cofibrant replacement (resp. fibrant replacement). It follows from
the axioms that a map $i$ is a cofibration iff it has the left lifting
property with respect to acyclic fibrations. Moreover $i$ is an
acyclic cofibration iff it has the left lifting property with
respect to fibrations. Therefore fibrations and weak equivalences
completely determine cofibrations. The dual properties hold for
fibrations.

\bigskip

Let $\cal C$ be a closed model category. We let $\cal W$ be the
class of all weak equivalences. Then there exists a category
${\cal H}o:={\cal C}[{\cal W}^{-1}]$ together with a functor
${\cal C}\to {\cal H}o$ which takes all morphisms from ${\cal W}$
to isomorphisms and which is universal with respect to this
property. Clearly  the category ${\cal H}o$  is determined
uniquely up to equivalence of categories. It has the following
description: objects of ${\cal H}o$ are the same as those of $\cal
C$, while morphisms are given by
$$\Hom_{{\cal H}o}(X,Y):=\Hom _{\cal C}(X^c,Y^f)/\sim ,$$ where
$\sim$ is an appropriate homotopy relation, which is  defined as
follows. Let $f,g:A\to B$ be two maps. Then $f\sim g$ if there
exists a map $h:IA \to B$ such that $f=h\circ i_1$ and $g=h\circ
i_2$. Here $IA$ and the maps $i_1,i_2:A\to IA$ satisfy the
following conditions: the canonical map   $(id,id):A\coprod A\to
A$ is a composite $A\coprod A{\buildrel (i_1,i_2)
\over\longrightarrow} IA \to A$, where the first map is a
cofibration and the second one is an acyclic fibration. It turns
out that this relation is an equivalence relation on $\Hom_{\cal
C}(A,B)$ if $A$ is cofibrant and $B$ is fibrant. Moreover it is
compatible with the composition law in $\cal C$ and the category
${\cal H}o$ is well defined.

\subsection{Cofibrantly generated model categories} Suppose $\cal C$ is a
category with all  colimits. Let $I$ be a class of maps in $\cal
C$. Following \cite{hov} we call a morphism {\it
$I$-injective} (resp. {\it $I$-projective})  if it has the right
(resp. left) lifting property with respect to every morphism in
$I$. The class of $I$-injective and $I$-projective morphisms are
denoted $I$-\emph{inj} and  $I$-\emph{proj} respectively. A morphism is
called an {\it $I$-cofibration} (resp. {\it $I$-fibration})  if it
has the left (resp. right) lifting property with respect to every
morphism in $I$-\emph{inj} (resp. $I$-\emph{proj}). The class of
$I$-cofibrations and $I$-fibrations are denoted $I$-\emph{cof} and
$I$-\emph{fib} respectively.  Assume now $I$ is a set of morphisms. A
morphism $f:A\to B$ is
 called a {\it relative $I$-cell complex} if there is an ordinal  $\lambda$ and
  a $\lambda$-sequence $$X_0\to X_1\to \cdots \to X_{\beta}\to \cdots, $$
$ \beta\leq \lambda$, with $A=X_0$ and $B=colim X_{\beta}$ such
that for all $\beta$ with $\beta + 1<\lambda$ there is a pushuot
diagram
$$\xymatrix{C_{\beta}\ar[r]^{g_{\beta}}\ar[d] & D_{\beta}\ar[d]\\
 X_{\beta}\ar[r]& X_{\beta +1}}$$
such that $g_{\beta}\in I$. The class of relative $I$-cell
complexes is denoted $I$-\emph{cell}. An object $A$ is called an {\it
$I$-cell complex} if $0\to A$ is a relative $I$-cell complex.

We will say that an object $A$ is {\it small}
relative to a class of morphisms $I$ if there
exists a cardinal $\kappa$ such that for each $\kappa$-filtered
ordinal  $\lambda$ and a $\lambda$-sequence
$X_0\to X_1\to \cdots \to X_{\beta} \to \cdots$ one has
$$colim \Hom_{\cal C}(A, X_{\beta})\cong \Hom_{\cal C}(A, colim X_{\beta}).$$
If $A$ is small with respect of $\cal C$ then $A$ is called {\it
small}. The following result is well-known (see for example
Theorem 2.1.19 of \cite{hov}).

 \begin{Pro}\label{Ho} Suppose $\cal C$ is a category with all  colimits and
 limits. Suppose $W$ is a subcategory of $\cal C$ and $I$ and $J$ are two sets
  of morphisms of $\cal C$ such that the following conditions hold:
\begin{itemize}
\item[(i)]  The subcategory $W$ is closed under retracts and
satisfies the CM2 axiom.

\item[(ii)] The domains of $I$ (resp. $J$) are small relative to
 $I$-\emph{cell} (resp. $J$-\emph{cell}).

\item[(iii)] $J$-\emph{cell} $\subset W\bigcap I$-\emph{cof}.

\item[(iv)] $I$-\emph{inj} $= W\bigcap J$-\emph{inj}.
\end{itemize}
Then there is a close model category structure on
$\cal C$ with $W$ as the subcategory of weak equivalences,
$I$-\emph{cof} as the class of cofibrations, $J$-\emph{inj} as the class  of
fibrations. Moreover $I$-\emph{inj} is the class of acyclic fibrations
and $J$-\emph{cof} is the class of acyclic cofibrations.
\end{Pro}

 The closed model categories obtained in this way are called {\it cofibrantly
 generated model categories}.

\subsection{Chain algebras} We fix a commutative ring $\K$ and all
algebras in what follows in this section are $\K$-algebras. Let us
recall that a chain algebra is a graded algebra
$A=\bigoplus_{n\geq 0} A_n$ equipped with a differential $d:A_n\to
A_{n-1}$ satisfying the Leibniz identity:
$$d(xy)=d(x)y+(-1)^{n}xd(y),   \ \ x\in A_n, y\in A_m.$$
Let ${\bf DGA}$ be the category of chain algebras.
%A morphism is called {\it weak equivalence} if it induces
%isomorphism on homology.

\begin{The}
 Define a morphism  of chain algebras to be
\begin{itemize}\label{CMCDGA}
\item[(i)] a weak equivalence if it induces isomorphism in homology
\item[(ii)] a fibration if it is a surjection in positive dimensions
\item[(iii)] a cofibration if it has the left lifting property
with respect to all maps which are fibrations and weak
equivalences
\end{itemize}
Then with these choices $\bf DGA$ is a cofibrantly generated
closed model category.
\end{The}

 To prove the theorem, we first introduce two classes of chain
 algebras. They play the role of discs and spheres. For
$n\geq 1$ we let $D(n)$ be the following chain algebra. As graded
algebra it is freely generated by elements $x$ and $dx$ of degree
$n$ and $n-1$ respectively. The boundary map assigns $dx$ to $x$.
For $n=0$ we let $D(0)$ be the algebra freely generated by an
element $x$ of degree $0$ (of course $d(x)=0$ in this case).
Moreover we define $S(n)$ to be the trivial algebra $\K$ if $n=-1$
and the algebra freely generated by an element $y$ of degree $n$
with zero boundary $d(y)=0$ provided $n\geq 0$. Then for all
$n\geq 0$ we have a canonical homomorphism $S(n-1) \to D(n)$ which
takes the generator $y$ to $dx$. We let $\coprod$ denote the
coproduct in $\bf DGA$. One has the following isomorphism of chain
complexes
\begin{equation}\label{copr}
A_*\coprod D(n)\cong A_*\oplus (A_*\t C_*\t A_*)\oplus (A_*\t
C_*\t A_*\t C_*\t A_*)\oplus \cdots.
\end{equation}
Here $C_*$ is a chain complex, which is zero in all dimensions
except for dimensions $n$ and $n-1$, where it is $\K$ and the unique
nontrivial boundary map is the identity. Therefore the inclusion
$A_*\to A_*\coprod D(n)$ is a weak equivalence, provided $n>0$.

One observes that for any chain algebra $A_*$ one has the
isomorphisms
\begin{equation}\label{Dn}\Hom_{\bf DGA}(D(n),A_*)\cong A_n,
\end{equation}
\begin{equation}\label{Sn} \Hom_{\bf DGA}(S(n),A_*)\cong \Ker(d:A_n\to A_{n-1}).
\end{equation}
We let $W$ be the subcategory of all weak equivalences in $\bf
DGA$. Moreover we put $$J:=\{\K\to D(n)\}_{n\geq 1},$$
$$I:=J\bigcup \{S(n-1)\to D(n)\}_{n\geq 0}.$$
 Then the conditions i) and ii) of Proposition \ref{Ho} hold, because
 $D(n)$ and $S(n)$ are small thanks to isomorphisms  \ref{Dn} and \ref{Sn}.
 We will show that all conditions of Proposition \ref{Ho} hold as well.
 To this end we need some preparations.

%We shall call a morphism $f:X_*\to Y_*$ {\it fibration} if it
%belongs to $J$-$inj$.
Since $\K$ is the initial object in $\bf DGA$ a morphism $f:X_*\to
Y_*$ is in $J$-\emph{inj} iff for any diagram
$$
\xymatrix{& X_*\ar[d]^f \\
 D(n)\ar[r]_{g} &  Y_*}$$ there exists a morphism $h:D(n)\to X$ such that $f=gh$. Now
the isomorphism (\ref{Dn}) gives that  $f:X_*\to Y_*$ is in
$J$-\emph{inj} iff  $f_n$ is surjective for all $n>0$. Thus we proved
the following
\begin{Le}\label{fib} A map $f:X_*\to Y_*$ is a $J$-$inj$ iff it is fibration.
\end{Le}

\begin{Le}\label{iinj} Let  $f:X_*\to Y_*$ be a
 morphism in $\bf DGA$. Then the following conditions are equivalent:
\begin{itemize}
 \item[(i)] $f$ is $I$-injective
 \item[(ii)]  $f$ is fibration and weak equivalence
 \item[(iii)]  $f_n$ is surjective for all $n\geq 0$ and $\Ker \, f$ is
acyclic, that is $H_*(\Ker \, f)=0$.
\end{itemize}
\end{Le}

\noindent {\it Proof}. Lemma \ref{fib} and the homology exact
sequence show that  iii) $\Longrightarrow $ ii).  Thanks to the
isomorphism (\ref{Sn}) a morphism  $f$ lies in $I$-\emph{inj} iff $f$
is a fibration with the following property: for all $x\in X_{n-1}$
and $y\in Y_n$ with $dx=0$ and $fx=dy$ there exists $z\in X_n$
such that $dz=x$ and $fz=y$. If the last condition holds, then
$f_0$ is surjective and $\Ker \, f$ is acyclic. Thus by Lemma
\ref{fib} we have i) $\Longrightarrow $ iii). Conversely, assume
iii) holds. Suppose $x\in X_{n-1}$ and $y\in Y_n$ are given with
$dx=0$ and $fx=dy$. Then there is $u\in X_n$ such that $fu=y$.
Since $f(x-du)=0$ and $d(x-du)=0$ it follows that $x-du=dv$ for
some $v\in X_n$ and therefore $x=d(u+v)$ which shows that  iii)
$\Longrightarrow $ i). To show  ii) $\Longrightarrow $ iii) it
suffices to show that $X_0\to Y_0$ is surjection. But this follows
from the commutative diagram
\bigskip
$$
\xymatrix{& \cdots \ar[r] & X_1\ar[r]\ar[d] &X_0\ar[r]\ar[d]
&H_0(X_*)\ar[r]\ar[d]_{\cong} &0\\
&\cdots \ar[r] & Y_1\ar[r]\ar[d] &Y_0\ar[r] &H_0(Y_*)\ar[r] &0\\
&&0&}
$$\rdg

\begin{Co}\label{ijinj} We have
 $I$-\emph{inj}$=W\bigcap J$-\emph{inj}.
\end{Co}

Recall that a morphism $f:X_*\to Y_*$ belongs to $I$-\emph{cof} if it
has the left lifting property with respect to all maps from
$I$-\emph{inj}. Thanks to Lemma \ref{iinj} this happens iff $f$ is a
cofibration.

A chain algebra $A_*$ is called $\K$-{\it projective} if each $A_n$ is
projective as a $\K$-module. A chain algebra $A_*$ is called {\it quasi-free} if
its underlying algebra is free. Let us recall that a graded algebra is free if it
is isomorphic to the tensor algebra $T(V)$ of a graded $\K$-module $V_*$, which
is free as a $\K$-module. A map of chain algebras $f:A_*\to B_*$ is called
{\it quasi-free} if $B_*$ as a graded algebra is a coproduct $A_*\coprod X_*$
where $X_*$ is a free algebra.

\begin{Le}\label{quasifree}{Quasi-free maps are  cofibrations.}\end{Le}
 {\it Proof}. Let
$$\xymatrix{A_*\ar[r]^{g}\ar[d]^i& X_*\ar[d]_p \\
B_*\ar[r]_h & Y_* }$$ be a commutative diagram of chain algebras,
in which $i$ is quasi-free and $p$ is an acyclic fibration. We
have to prove that there is a chain map $f:B_*\to X_*$  such that
$g=fi$ and $h=pf$. By assumptions we have an isomorphism of
algebras $B_*\cong A_*\coprod C_*$, where $C_*$ is a free algebra.
Let $E$ be the set of free generators of $C_*$. Then $E$ is the
union of subsets $E_n$ of degree $n$ elements, $n\geq 0$. In order
to define $f$ one needs to specify  elements $f(x)$ for $x\in
E_n$, $n\geq 0$ with two properties
\begin{itemize}
\item[a)]  $\partial f(x)= f(\partial x)$,
\item[b)] $pf(x)= h(x)$.
\end{itemize}
We will work by induction on $n$. First consider the case $n=0$.
Since $p$ is surjective, there exists $f(x)\in X_0$ such that
$pf(x)=h(x)$. Consider now the case $n>0$. Suppose  for all $m<n$
we already defined $f(x)$ for all $x\in E_m$  such that a) and b)
holds for all $x\in E_j$, $1\leq j<n$. Take now $x\in E_n$. Since
$p$ is surjective we can choose an element ${\bar f}(x)\in X_n$
such that $p{\bar f}(x)= g(x)$. Since $\partial x$ lies in the
subalgebra generated by $A_*$ and $E_j$, $j<n$, the element
$f(\partial x)$ is already defined. Set $z=
\partial ({\bar f}(x))-f(\partial (x))$. Then $\partial (z)= 0$
and $p(z)=0$. Therefore $z=\partial (u)$ for some element $u\in
\Ker (p)$. Now we put $f(x)= {\bar f}(x)-u$. It is clear that
$f(x)$ satisfies properties a)  and b). Thus induction step is
finished and hence the lemma.

\rdg

\Cor{ The canonical maps $k\to D(n)$, $k\to S(n)$ and $S(n-1)\to
D(n)$ are cofibrations.}

%\begin{The} The category $\bf DGA$ equipped with weak equivalences,
%cofibrations and fibrations is a cofibrantly generated model
%category. \end{The}
%: Define a morphism $f:A_*\to B_*$ to be
%\S
%{\rm (i)} a weak equivalence
%if $f$ induces isomorphism in homology%
%\S
%{\rm (ii)} a fibration if it is surjection in positive dimensions
%\S
%{\rm (iii)} a cofibration if $f$ has the left lifting property with respect to
%all maps which are fibrations and weak equivalences
%
%Then with these choices $\bf DGA$ is a closed model category.}
 {\bf Proof of Theorem \ref{CMCDGA}.} As was mentioned already the conditions i) and ii)
of Proposition \ref{Ho} hold. By Corollary \ref{ijinj} the condition iv) holds
as well. Thus we have to show that $J$-\emph{cell}$\subset W\bigcap I$-\emph{cof}. We have
$J$-\emph{cell}$\subset I$-\emph{cof} because $J\subset I$. Since the domain of all maps
from $J$ is $\K$, which is an initial object, we see that all pushouts in the
definition of a relative $J$-cell complex are coproducts with $D(n)$ for some
$n>0$. It follows that all such morphisms are weak equivalences and quasi-free
maps and the result follows from Lemma \ref{quasifree}.

Let us note that a similar theorem
for cochain algebras was proved  by Jardine \cite{jar}. Moreover
our argument is merely a variant of the one given there (compare also
with \cite{bg}).

\subsection{Truncated chain algebras}\label{cmctrun} Let us fix a natural number
$m\geq 1$. We let ${\bf DGA}_m$ be the full subcategory of $\bf
DGA$ which consists of objects $X_*$ such that $X_i=0$ for all
$i>m$.

For any chain complex $(X_*,d)$ we let $\tam(X_*)$ be the
following chain complex:
$$(\tam(X_*))_i=X_i, \ {\rm if } \ i<m$$
$$(\tam(X_*))_m=X_m/d(X_{m+1})  \ \ $$
$$(\tam(X_*))_i=0, \ {\rm if } \ i>m$$
The quotient map $X_*\to \tam(X_*)$ is a chain map. Moreover
$H_i(\tam(X_*))\cong H_i(X_*)$ if $i\leq m$ and $H_i(\tam(X_*))=0$
provided $i>m$. It is also clear that, if $X_*$ is a chain
algebra, then there is a unique chain algebra structure on
$\tam(X_*)$ such that the quotient map $X_*\to \tam(X_*)$ is a
chain algebra homomorphism. Thus $$\tam:{\bf DGA}\to {\bf DGA}_m$$
is a well-defined functor, which is the left adjoint  to the
inclusion functor ${\bf DGA}_m\subset {\bf DGA}$.

\begin{The}\label{CMCDGAM} Define a map  in ${\bf DGA}_m$ to be a weak
equivalence (resp. fibration) if  it is a weak equivalence (resp.
fibration) in $\bf DGA$. Define a map  in ${\bf DGA}_m$ to be a
cofibration if it has left lifting property with respect to all acyclic
fibrations. Then this defines a cofibrantly generated model
structure on ${\bf DGA}_m$.
\end{The}

{\it Proof}. We introduce two classes of morphisms in ${\bf
DGA}_m$:
$$J_m:=\{\K\to \tam D(n)\}_{ n\geq 1 },$$
$$I_m:=J_m\bigcup \{\tam S(n-1)\to \tam D(n)\}_{n\geq 0}.$$
We have to show that all assertions of  Proposition \ref{Ho}
hold. Conditions i) and ii) are clear. Formal argument with
adjoint functors shows that a morphism
 $f:X_*\to Y_*$ in  ${\bf DGA}_m$ considered as a morphism of $\bf DGA$ lies in
  $J$-\emph{inj} (resp. $I$-\emph{inj})
 iff it is in $J_m$-\emph{inj} (resp. $I_m$-\emph{inj}). Therefore $f$ is a
 fibration (resp. acyclic fibration) iff it is in $J_m$-\emph{inj} (resp.
 $I_m$-\emph{inj}) and the condition iv) holds.
 We also have $J_m$-\emph{cell}$\subset  I_m$-\emph{cof} because $J_m\subset I_m$.
 Thus it remains to show that $J_m$-\emph{cell}$\subset W$. Comparing the
 definitions we see that any morphism from $J_m$-\emph{cell} can be
 written as $\tam(g)$, where $g\in J$-\emph{cell}. In particular $g\in W$.
  Since $\tam$ preserves weak equivalences we are done.

\subsection{A closed model category structure on crossed bimodules}\label{cmcxmod} Of the
special interest is the case, when $m=1$. In this case Theorem \ref{CMCDGAM}
gives the closed model category structure on the category $\sf Xmod$ of crossed
bimodules. A map of crossed bimodules
$$\xymatrix{C_1\ar[r]^\d\ar[d]^f&C_0\ar[d]^g\\
C_1' \ar[r]^{\d'}&C_0'}
$$
is a fibration if $f$ is a surjective homomorphism. Moreover, $(f,g): \d\to
\d'$ is a weak equivalence if induced maps $\ker(\d)\to \ker(\d')$ , ${\sf
Coker}(\d)\to \ {\sf Coker}(\d')$ are isomorphisms. It follows that if $(f,g)$
is an acyclic fibration, then $g$ is a surjection and the induced map
$\ker(f)\to \ker(g)$ is an isomorphism, in other words one has the following
commutative diagram with exact rows and columns:
$$\xymatrix@=1em{&& 0\ar[d]& 0\ar[d]&&\\
&&\ker(f)\ar[d]\ar[r]^{\cong}&\ker(g)\ar[d]&&\\
0\ar[r]&\ker(\d)\ar[d]^\cong
\ar[r]&C_1\ar[d]^f\ar[r]^\d&C_0\ar[r]\ar[d]^g&{\sf Coker}(\d)\ar[r]\ar[d]^\cong & 0\\
0\ar[r]&\ker(\d')
\ar[r]&C_1'\ar[r]^{\d'}\ar[d]&C_0\ar[r]\ar[d]&{\sf Coker}(\d')\ar[r]& 0\\
&&0&0&&}
$$
Thus we proved the following
\begin{Le} If $(f,g):\d\to \d'$ is an acyclic fibration in $\sf Xmod$, then $g$ is
surjective and
$$\xymatrix{
 C_1\ar[d]^f\ar[r]^\d&C_0\ar[d]^g\\
C_1'\ar[r]^{\d'}&C_0'}
$$
is a pullback diagram.
\end{Le}
  \begin{Le}
  A crossed bimodule $\delta :R_1\to R_0$ is
  a cofibrant objects in $\sf Xmod$ provided $R_0$ is a free algebra.
\end{Le}

{\it Proof}.  Let $(f,g):\d\to \d'$ be an acyclic fibration of crossed
bimodules and let $(a',b'):\delta\to \d'$ be a morphism of crossed bimodules. We
have to lift it to a morphism $(a,b):\delta\to \d$. Since $g$ is a surjective
homomorphism of $\K$-algebras and $R_0$ is a free $\K$-algebra, we can lift
$b'$ to a homomorphism $b:R_0\to C_0$ of $\K$-algebras. Then we have the
following commutative diagram
$$
\xymatrix@=.2em{R_1\ar[dddddr]_{a'}\ar[drrrrr]^{b\partial }\\
& C_1\ar[dddd]^f\ar[rrrr]_\d&&&&C_0\ar[dddd]^g\\
\\
\\
\\
&C_1'\ar[rrrr]^{\d'}&&&&C_0'}
$$
and we have the unique homomorphism $a:R\to C_1$ which fits in the diagram. It
is now clear that $(a,b)$ is an expected lifting.


\begin{thebibliography}{33}

\bibitem{arabia} A. Arabia, Rel\'evements des alg\'ebres lisses et de leurs
morphismes. Comment. Math. Helv. 76(2001) 607-639.

\bibitem{BW} H.-J. Baues and G.Wirsching, Cohomology of small categories.
J. Pure Appl. Algebra 38 (1985), no. 2-3, 187--211.

\bibitem{B1} H.-J. Baues. Combinatorial homotopy
and $4$-dimensional complexes. With a preface by Ronald Brown. de Gruyter
Expositions in Mathematics, 2. Walter de Gruyter., Berlin, 1991. xxviii+380 pp.

\bibitem{baumin} H.-J. Baues and E.C. Minian. Crossed extensions of algebras
and Hochschild cohomology. Homology, homotopy and applications,
4(2002) 63-82.

\bibitem{BJP} H.-J. Baues, M. Jibladze and T. Pirashvili.
Strengthening  track theories. arXiv entry \texttt{math.CT/0307185}.

\bibitem{square} H.-J. Baues, M. Hartl  and T.Pirashvili. Quadratic categories and square rings. J.
  Pure Appl. Algebra 122 (1997) 1-40.


\bibitem{bo} M. B\"okstedt.  The topological Hochschild homology of
$\Z$ and $\Z/p$. Unpublished manuscript.

\bibitem{bg} A.K. Bousfield and V. K. A. M. Gugenheim,
On ${\rm PL}$ de Rham theory and rational homotopy type. Mem.
Amer. Math. Soc. 8 (1976), no. 179, ix+94 pp.

\bibitem{brun} M. Brun. Topological Hochschild homology of $\Z/p^n\Z$ J. Pure
Appl. Algebra, 148(2000),29-76.

\bibitem{CE} {\sc H. Cartan} and {\sc S. Eilenberg}. Homological algebra. Princeton University Press, Princeton, N. J., 1956. xv+390 pp.

\bibitem{delu} P. Dedecker and A. S.-T. Lue. A nonabelian two-dimensional
cohomology for associative algebras. Bull. AMS 72(1966), 1044-1050.

\bibitem{dum} B. I. Dundas and R. McCarthy.
  Stable $K$-theory and topological Hochschild homology.  Ann. of Math.
  (2)  140  (1994),  no. 3, 685--701.  Erratum: 142  (1995),  no. 2, 425--426.

\bibitem{duskin} J. Duskin. Simplicial methods and the interpretation of
"triple" cohomology. Mem. AMS 163(1975).

\bibitem{dwyerspalinski}  W.G. Dwyer and J.Spalinski,  Homotopy theories
and model categories. Handbook of algebraic topology, 73--126,
North-Holland, Amsterdam, 1995.

\bibitem{EM}   S. Eilenberg and   S. MacLane, Homology theory
for multiplicative systems. Trans. AMS 71  (1951), 294-330.


 \bibitem{fls} V.Franjou, J. Lannes and L. Schwartz.
Autour de la cohomologie de MacLane des corps finis,  Inven. Math.
115, (1994), 513-538.


\bibitem{FP} V.Franjou and T. Pirashvili.
On the Mac Lane cohomology for the ring of integers. Topology 37
(1998), 109--114.

\bibitem{gerstenhaber} M. Gerstenhaber. The cohomology structure of an
associative ring. Ann. of Math. 78(1963), 267-288.

\bibitem{hov} M. Hovey. Model categories. Mathematical Survays and
Monographs, 63, AMS, Providance, RI,1999,xii+209pp.

\bibitem{hueb} J. Huebschmann. Crossed $n$-fold extensions of
groups and cohomology. Comment. Math. Helvetici. 55(1980),302-314.

\bibitem{jar} J.F. Jardine. A closed model structure for
differential graded algebras. Cyclic cohomology and noncommutative
geometry (Waterloo, ON, 1995), 55--58, Fields Inst. Commun., 17,
Amer. Math. Soc., Providence, RI, 1997.

\bibitem{JP1} M.Jibladze and T. Pirashvili.
Some linear extensions of a category of finitely generated free
modules. Soobshch. Akad. Nauk Gruzin. SSR 123 (1986), 481--484.

\bibitem{JP2} M.Jibladze and T. Pirashvili. Cohomology of
algebraic theories. J. Algebra 137 (1991), no. 2, 253--296.

\bibitem{jll} J.-L. Loday. Spaces with finitely many nontrivial
homotopy groups.  J. Pure Appl. Algebra  24  (1982), 179--202.

\bibitem {jllst} {\sc J. L. Loday},  {\it Cohomologie et
groupe de Steinberg relatifs.} J.  Alg. {\bf 54} (1978), 178-202.

\bibitem{HC} J.-L. Loday. Cyclic homology. Second edition.
Grundlehren der Mathematischen Wissenschaften. 301.
Springer-Verlag, Berlin, 1998. xx+513 pp.

\bibitem{mac} S. MacLane. {\it Homologie des anneaux et des modules}, Coll.
 topologie algebrique, Louvan (1956), 55-80.

\bibitem{macobs} S. MacLane. Extensions and obstructions for
rings. Ill. J. of Math. 2(1958), 316-345.

\bibitem{homology} S. MacLane. {\it Homology},
Springer-Verlag Berlin, Heidelberg, New York, 1975.

\bibitem{MW} S. MacLane and J.H.C. Whitehead. On the 3-type of a
complex. Proc. Nat. Acad. of Sci. USA, 36(1950), 41-48.


\bibitem{qu} D. G. Quillen. \emph{Homotopical algebra}. Springer
Lecture Notes in Mathematics, {\bf 43}, Springer-Verlag, 1967.

\bibitem{Q} {\sc D. Quillen}. On the (co-) homology of commutative rings.
Proc. Sympos. Pure Math., Vol. XVII, New York, (1968) 65--87.


\bibitem{P1} T. Pirashvili.
Models for the homotopy theory and cohomology of small categories. Soobshch.
Akad. Nauk Gruzin. SSR 129 (1988), 261--264.

 \bibitem{P2} T. Pirashvili. Cohomology of small categories
in homotopical algebra. $K$-theory and homological algebra (Tbilisi, 1987--88),
268--302, Lecture Notes in Math., 1437, Springer, Berlin, 1990.

\bibitem{P3} T. Pirashvili. Spectral sequence for Mac Lane homology.
J. Algebra 170 (1994), 422--428.

\bibitem{P4} T. Pirashvili.
On the topological Hochschild homology of $\Z/p^k\Z$. Comm.
Algebra 23 (1995), 1545--1549.

\bibitem{PW} T. Pirashvili and F.Waldhausen.
Mac Lane homology and topological Hochschild homology. J. Pure
Appl. Algebra 82 (1992), 81--98.

\bibitem{relshukla}T. Pirashvili. Algebra cohomology
over a commutative algebra revisited. math.arXive: KT/0309184.

\bibitem{sh} U.Shukla. Cohomologie des alg\'ebres associatives.  Ann. Sci. \'Ecole Norm. Sup. (3) 78 (961)
163--209.

\bibitem{yoneda} N. Yoneda. On $Ext$ and exact sequences. J. Fac.
Sci. Univ. Tokyo. Sect. I, 8,1960, 507-576(1960).

\bibitem{hanry} J. H. C. Whitehead. Combinatorial homotopy II.
Bull. AMS 55(1949), 453-496.

\end{thebibliography}
\end{document}